\documentclass[a4paper,preprint,10pt,3p]{elsarticle}

\RequirePackage{mathtools}  
\RequirePackage{amssymb}    
\RequirePackage{siunitx}    

\RequirePackage{tabularx}   
\RequirePackage{booktabs}   
\RequirePackage{longtable}  
\RequirePackage{multirow}   

\RequirePackage{graphicx}   
\RequirePackage{float}      
\RequirePackage[labelfont=bf,justification=centering,footnotesize]{caption} 
\RequirePackage{subcaption} 
\RequirePackage{pdfpages}   

\RequirePackage{xcolor}     
\RequirePackage{tikz}       
\RequirePackage{xspace}     
\RequirePackage{microtype}  

\RequirePackage{geometry}   
\RequirePackage{titlesec}   
\RequirePackage{titletoc}   
\RequirePackage{fancyhdr}   
\RequirePackage{enumitem}   
\RequirePackage{etoolbox}   
\RequirePackage{iftex}      
\RequirePackage{datetime}   
\RequirePackage{hyperref} 
\RequirePackage[noabbrev,nameinlink]{cleveref} 

\RequirePackage{listings} 

\usepackage{amsmath}

\usepackage{algorithm}
\usepackage[noend]{algpseudocode}
\usepackage{comment}
\usepackage{placeins}

\DeclareMathOperator*{\argmin}{arg\,min}

\newcommand{\bm}[1]{\boldsymbol{#1}}

\usepackage{amsthm}
\usepackage{xpatch}

\xpatchcmd{\proof}{\itshape}{\prooflabelfont}{}{}
\newcommand{\prooflabelfont}{\bfseries}

\newtheorem{definition}{Definition}

\usepackage{array}
\newcolumntype{L}[1]{>{\raggedright\let\newline\\\arraybackslash\hspace{0pt}}m{#1}}

\usepackage[pagewise]{lineno} 

\newdefinition{rmk}{Remark}

\begin{document}
\begin{frontmatter}

\title{Entropy-Stable Model Reduction of One-Dimensional Hyperbolic Systems using Rational Quadratic Manifolds}
\author[1,2]{R.B. Klein\corref{cor1}}
\ead{rbk@cwi.nl}
\author[1,3]{B. Sanderse}
\ead{B.Sanderse@cwi.nl}
\author[2]{P. Costa}
\ead{P.SimoesCosta@tudelft.nl}
\author[2]{R. Pecnik}
\ead{R.Pecnik@tudelft.nl}
\author[2]{R.A.W.M. Henkes}
\ead{R.A.W.M.Henkes@tudelft.nl}

\cortext[cor1]{Corresponding author}


\address[1]{Centrum Wiskunde \& Informatica, Science Park 123, Amsterdam,   The Netherlands}
\address[2]{Delft University of Technology, Process and Energy, Leeghwaterstraat 39, Delft, The Netherlands}
\address[3]{Eindhoven University of Technology, Department of Mathematics and Computer Science, PO Box 513, Eindhoven, The Netherlands}

\begin{abstract}
In this work we propose a novel method to ensure important entropy inequalities are satisfied semi-discretely when constructing reduced order models (ROMs) on nonlinear reduced manifolds. We are in particular interested in ROMs of systems of nonlinear hyperbolic conservation laws. The so-called entropy stability property endows the semi-discrete ROMs with physically admissible behaviour. The method generalizes earlier results on entropy-stable ROMs constructed on linear spaces. The ROM works by evaluating the projected system on a well-chosen approximation of the state that ensures entropy stability. To ensure accuracy of the ROM after this approximation we locally enrich the tangent space of the reduced manifold with important quantities. Using numerical experiments on some well-known equations (the inviscid Burgers equation, shallow water equations and compressible Euler equations) we show the improved structure-preserving properties of our ROM compared to standard approaches and that our approximations have minimal impact on the accuracy of the ROM. We additionally generalize the recently proposed polynomial reduced manifolds to rational polynomial manifolds and show that this leads to an increase in accuracy for our experiments.
\end{abstract}


\end{frontmatter}

\section{Introduction}
Conservation laws are nearly universally present in any branch of physics and engineering e.g.\ fluid dynamics, structural mechanics, plasma physics and climate sciences; they express the conservation of some physical quantity of interest. Often such conservation laws are described by hyperbolic equations or systems thereof \cite{dafermoshyperbolic}. Physicists and engineers are increasingly reaching to simulation tools for approximate solutions. With the increase in computational power of recent decades, very large-scale problems have indeed been solved to a satisfactory accuracy. Nonetheless, some applications of major engineering interest like those of a multi-query (e.g.\ design optimization \cite{rozzabasic}, uncertainty quantification \cite{deliaquantification}) or real-time nature (e.g.\ model predictive control \cite{noackreduced}, digital twin technology \cite{hartmannmodel}) remain out of question for many large scale systems. A way around these computational issues has been offered by reduced order models (ROMs), which are low-dimensional surrogates of high-fidelity models of interest, often referred to as full order models (FOM) in the ROM community. ROMs rely on the offline-online decomposition paradigm \cite{lassilamodel} for their efficiency, they are trained in an expensive offline phase and subsequently evaluated at very low computational cost in new situations in the online phase. Particularly, the low dimensionality of the ROM allows for fast and cheap evaluation of its solution.

A popular class of ROMs are projection-based ROMs (pROMs) \cite{bennersurvey, rozzabasic}, which have traditionally been constructed by projecting equations of interest on well-chosen linear subspaces. These subspaces are often found in a data-driven manner using the proper orthogonal decomposition (POD) \cite{turbulenceholmes, sirovichturbulence} or greedy methods \cite{quarteronireduced} and the projections are carried out using a Galerkin \cite{sandersenonlinearly, kleinenergy, rowleymodel, kalashnikovastable, arunajatesanstable, baronereduced} or Petrov-Galerkin \cite{carlbergefficient, grimbergmesh, carlbergmodel} approach. However, the success of applications of linear subspace-based pROMs has been limited in the field of hyperbolic equations. This is a result of an almost inherently slow decay of the so-called Kolmogorov $n$-width (KnW) of the solution manifolds $\mathcal{M}_u$ (i.e.\ the set of all solution trajectories for a range of parameters and initial conditions of interest) of these systems\footnote{$(\mathcal{X},||\cdot||_{\mathcal{X}})$ denoting a Banach space of interest containing analytical PDE or approximate simulation solutions}:
\begin{equation}
    d_r(\mathcal{M}_u) = \inf_{\mathcal{V}\subset\mathcal{X}; \text{dim}(\mathcal{V}) = r} \sup_{\bm{u}\in\mathcal{M}_u} \inf_{\bm{v}\in\mathcal{V}} ||\bm{u} - \bm{v}||_{\mathcal{X}},
    \label{eq:knW}
\end{equation}
which measures the worst error that can be incurred when optimally approximating $\mathcal{M}_u$ with an optimally chosen $r$-dimensional linear subspace. The KnW decay of many hyperbolic systems is slow because their solution trajectories are often not contained in low-dimensional linear subspaces due to characteristically having moving features as part of their solution. This has been shown analytically for some systems \cite{pinkusn, cohenkolmogorov, greifdecay, melenkn, bachmayrkolmogorov} and empirically for many others. In recent years a range of possible solutions have been proposed, falling in roughly four categories. First, there have been adaptive approaches that use linear subspaces that are changed during the online phase to be better suited to new conditions \cite{peherstorferonline, peherstorfermodel, pagliantinifully, hesthavenrank, carlbergadaptive, etteronline, zahrefficient, gaoadaptive}; second, domain decomposition approaches that localize ROM construction in time, space or parameter space \cite{washabaughnonlinear, chaturantabuttemporal, peherstorferlocalized, ijzermansignal, borggaardinterval, ahmedbreaking, grimbergmesh, copelandreduced}; third, Lagrangian, registration and/or optimal transport based approaches that track moving features improving linear data compressibility \cite{ohlbergernonlinear, gerbeauapproximated, iolloadvection, reissshifted, mojganilagrangian, taddeiregistration, ehrlachernonlinear, rimmanifold, mojganilow, khamlichoptimal, vanheyningenadaptive, blickhanregistration, mirhoseinimodel}; fourth, constructing ROMs on nonlinear spaces (manifolds). Since it is not dependent on user-defined localization or adaptation strategies and due to its high expressiveness we will be interested in the latter category. In model reduction on manifolds the linear reduced spaces of classical ROMs are replaced by nonlinear spaces. Given sufficient expressiveness the nonlinear reduced spaces can potentially approximate the solution manifolds of hyperbolic systems. These manifolds are typically constructed from data. In \cite{romorexplicable, romornonlinear, leemodel, buchfinksymplectic, tencerenabling, cocolahyper, kimfast, diazfast} autoencoder neural networks are used to construct manifold ROMs. Another popular nonlinear manifold construction approach are polynomial manifolds \cite{geelenoperator, buchfinkapproximation, sharmasymplectic, barnettquadratic, jainquadratic}. Finally, there has been an increasing interest in a more physics-based approach, where ROMs are constructed on the invariant manifolds of physical systems \cite{axasmodel, buzaspectral, ottolearning, buzausing}.

Although, to the authors' knowledge, there have not been any studies of manifold ROMs on large scale under-resolved and shock-dominated cases within fluid dynamical applications, it has long been known that their linear counterparts can suffer from stability issues for such problems \cite{sandersenonlinearly, kleinenergy, fickstabilized, balajewiczstabilization, amsallemstabilization, kalashnikovastable, arunajatesanstable, baronereduced}. It is quite reasonable to assume this will also be the case for manifold ROMs. A promising approach to stabilization of such simulations is the concept of entropy stability \cite{tadmorentropy,leveque}, which has been widely used in the context of obtaining stable FOMs. A numerical method is entropy stable if it dissipates a convex functional associated to a conservation law, referred to as entropy, given suitable boundary conditions. Entropy stability endows numerical methods with physically admissible behaviour, which for fluid dynamical applications manifests itself in satisfaction of the second law of thermodynamics. Furthermore, it generalizes the $L^2$-stability properties of linear hyperbolic systems to fully nonlinear systems \cite{tadmorentropy}. Additionally, stability in $L^p$ spaces can be shown formally \cite{svardweak, fjordholmarbitrarily, kruzkovfirst}. 

\begin{figure}
    \centering
    \includegraphics[width = 0.7\textwidth]{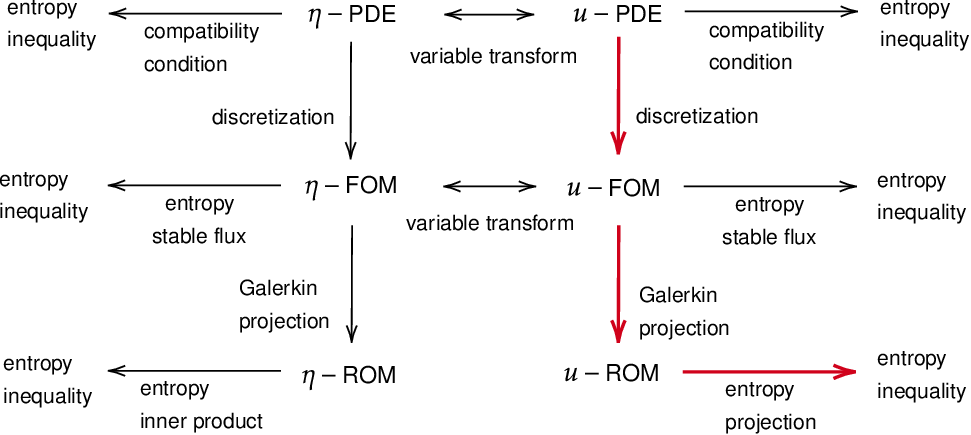}
    \caption{Flowchart of possible entropy stable ROM approaches, our method is indicated in red, $u$ are the conserved variables and $\eta$ are the alternate variables used in \cite{kalashnikovastable}. FOM: full order model, ROM: reduced order model.}
    \label{fig:flowchart}
\end{figure}

However, in the projection step to construct pROMs an entropy-stable numerical method generally loses its stability property. In recent years, some work has been carried out in the ROM community to preserve the entropy stability property. In \autoref{fig:flowchart} several possible approaches have been visualized. In \cite{chanentropy} the ROM is evaluated at a corrected state that ensures entropy stability is maintained (also known as `entropy projection'). Inspired by the approach taken in the classical finite element work in \cite{hughesnew}, \cite{kalashnikovastable} writes the conservation laws in symmetric form using an alternative set of variables and projects the continuous equations leading to correct entropy estimates. In \cite{parishimpact} the Hilbert spaces in which the projections are carried out are defined according to physical arguments. We also note the works in \cite{arunajatesanstable, baronereduced, rowleymodel}. All these approaches suffer a major drawback -- they are built on linear reduced spaces. As a result, relatively high-dimensional reduced spaces are required to model many physical systems of interest, which comes at the cost of computational efficiency of the ROM.

Our main contribution is to generalize these entropy-stable approaches to nonlinear reduced spaces. This allows for lower-dimensional reduced spaces and thus potentially more efficient ROMs. In particular, we will be interested in generalizing the work in \cite{chanentropy}. This approach is colored red in \autoref{fig:flowchart}. Our main argument for not choosing \cite{kalashnikovastable, parishimpact} stems from the argument given e.g.\ in \cite{carpenterentropy}, namely that formulations in alternative variables are not consistent with the famous Lax-Wendroff theorem \cite{laxhyperbolic} and can thus yield wrong shock solutions. To our knowledge, our method is the first manifold ROM that is provably entropy stable. At the same time, we note that preservation of other mathematical structures on reduced manifolds has been successfully achieved in the past: symplectic \cite{sharmasymplectic, buchfinksymplectic, brantnerstructurepreserving}; metriplectic \cite{leemachine}; conservative \cite{leedeep}. A great overview of recent structure-preserving model reduction contributions in general is given in \cite{grubercanonical}.

A second contribution of this work is the development of a novel generalization of polynomial manifolds \cite{barnettquadratic, buchfinkapproximation} to rational polynomial manifolds. While these polynomial manifolds have shown successes in certain applications, in our experience they are not sufficiently accurate for shock-dominated problems. Rational polynomials are better at capturing discontinuities, as we will show in this work.

This article is organized as follows. In \autoref{sec:prelims} we introduce the theory of nonlinear hyperbolic conservation laws and entropy analysis, we introduce a baseline entropy stable FOM \cite{fjordholmarbitrarily}, and an entropy-stable linear ROM as proposed in \cite{chanentropy}. In \autoref{sec:rom} we introduce our main contribution, the novel entropy stable nonlinear manifold ROM. In \autoref{sec:rational} we discuss our second contribution, being rational polynomial manifolds. In \autoref{sec:experiments} we show the effectiveness of our approach using several numerical experiments that are based on a range of well-known conservation laws from fluid dynamics. We conclude our work in \autoref{sec:conclusion}.


\section{Preliminaries: Entropy inequality for conservation laws, entropy-stable FOM, linear ROM} \label{sec:prelims}
\subsection{Introducing the entropy inequality}
We give a short introduction to the concept of entropy, some related concepts used in its analysis and its role in the theory of nonlinear conservation laws. We consider conservation laws in one spatial dimension that can be written as partial differential equations (PDE) of the form:
\begin{equation}
    \frac{\partial \bm{u}}{\partial t} + \frac{\partial \bm{f}(\bm{u})}{\partial x} = 0,
    \label{eq:conslaw}
\end{equation}
where $\Omega$ is a spatial domain and $[0,T]$ is a temporal domain with $T > 0$. Furthermore $\bm{u} : \Omega \times (0,T] \rightarrow \mathbb{R}^n$ is the solution function, $\bm{f} : \mathbb{R}^n \rightarrow \mathbb{R}^n$ is the nonlinear flux function, $n \in \mathbb{N}$ is the number of conserved quantities and $x\in \Omega$ and $t \in [0,T]$ are the spatial coordinate and time, respectively. To facilitate conservation statements and minimize the role of boundary conditions we will focus on periodic spatial domains $\Omega = \mathbb{T}([a,b])$ with $b > a$ and $b,a \in \mathbb{R}$ in this research, here $\mathbb{T}$ is the torus. The equations are complemented by a set of initial conditions $\bm{u}_0 : \Omega \rightarrow \mathbb{R}^n$ so that $\bm{u}_0(x) = \bm{u}(x,0)$. The conservation law \eqref{eq:conslaw} is referred to as hyperbolic when the Jacobian matrix $\frac{\partial \bm{f}}{\partial \bm{u}}$ is diagonalizable with real eigenvalues for all physically relevant $\bm{u}$ \cite{leveque}. In many cases, the solutions of physically relevant hyperbolic conservation laws also satisfy \textit{additional} conservation laws of the form:
\begin{equation}
    \frac{\partial s(\bm{u})}{\partial t} + \frac{\partial \mathcal{F}(\bm{u})}{\partial x} = 0,
    \label{eq:entconslaw}
\end{equation}
where the function $s : \mathbb{R}^n \rightarrow \mathbb{R}$ is called the entropy function which is defined to be convex\footnote{A function $g : \mathbb{R}^n \rightarrow \mathbb{R}$ is convex if its Hessian, $\frac{\partial^2 g}{\partial \bm{u}^2}(\bm{u})$, is positive definite for all $\boldsymbol{u}$.} and $\mathcal{F} : \mathbb{R}^n \rightarrow \mathbb{R}$ is called the entropy flux. In particular, such an additional conservation law exists if the compatibility relation:
\begin{equation}
    \bm{\eta}(\bm{u})^T\frac{\partial \bm{f}}{\partial \bm{u}} = \frac{\partial \mathcal{F}}{\partial \bm{u}}^T,
    \label{eq:comprelation}
\end{equation}
is satisfied. Here, $\bm{\eta} : \mathbb{R}^n \rightarrow \mathbb{R}^n, \hbox{ } \bm{u} \mapsto \frac{\partial s}{\partial \bm{u}}(\bm{u})$ is the gradient of the entropy function $s$ with respect to $\bm{u}$. This mapping is injective due to the convexity of $s$ and hence can be inverted. A pair $(s,\mathcal{F})$ satisfying the compatibility relation \eqref{eq:comprelation} is called an entropy pair of \eqref{eq:conslaw}.

It is well-known that solutions $\bm{u}$ to \eqref{eq:conslaw} can develop discontinuities in finite time for smooth $\bm{u}_0$ \cite{evanspartial, leveque, leflochhyperbolic}. When this occurs the solution $\bm{u}$ is said to contain a shock or contact discontinuity depending on the behaviour of the discontinuity \cite{leveque}. In this case both formulation \eqref{eq:conslaw} and the manipulations to obtain \eqref{eq:entconslaw} are no longer valid. We must therefore consider \eqref{eq:conslaw} in a weak sense to retain a notion of solutions. This weak form of the conservation law is obtained by integrating against a space of smooth test functions $\bm{v} : \Omega \times [0,\infty ) \rightarrow \mathbb{R}^n$ with compact support i.e.\ $\bm{v} \in C_0^{\infty}(\Omega \times [0,\infty) )$ and transferring all derivatives to these test functions to obtain:
\begin{equation}
    \int_{0}^{\infty} \int_{\Omega} \bm{u}\cdot \frac{\partial \bm{v}}{\partial t} + f(\bm{u}) \cdot \frac{\partial \bm{v}}{\partial x} dx dt + \int_{\Omega} \bm{u}_0 \cdot \bm{v}(x,0) dx = 0,
    \label{eq:weakform}
\end{equation}
on periodic $\Omega$. Note that this expression is valid even for discontinuous $\bm{u}$ and that any smooth $\bm{u}$ satisfying the strong form \eqref{eq:conslaw} also satisfies this weak form \eqref{eq:weakform} \cite{dafermoshyperbolic, evanspartial}. However, weakening the notion of a solution like this comes at the cost that it does not necessarily yield unique solutions (examples of such cases may be found in \cite{evanspartial, leflochhyperbolic, dafermoshyperbolic, leveque}). Out of all weak solutions, i.e.\ solutions to \eqref{eq:weakform}, those of physical interest are the ones 
satisfying:
\begin{equation}
    \frac{\partial s(\bm{u})}{\partial t} + \frac{\partial \mathcal{F}(\bm{u})}{\partial x} \leq 0,
    \label{eq:ineqentconslaw}
\end{equation}
in the sense of distributions, with equality for smooth $\bm{u}$ following \eqref{eq:entconslaw} and inequality for solutions containing shocks. This inequality arises from considering limits of regularized conservation laws, more details can be found in \cite{leveque}. For scalar conservation laws Kruzkov \cite{kruzkovfirst} established that weak solutions satisfying \eqref{eq:ineqentconslaw} are unique, but for systems uniqueness is not yet completely established \cite{bianchinivanishing}. We can define the total entropy functional as:
\begin{equation}
    \mathcal{S}[\bm{u}] := \int_{\Omega} s(\bm{u}) dx.
    \label{eq:entropycont}
\end{equation}
Defining appropriate sequences of test functions and taking limits \cite{levequenumerical, dafermoshyperbolic}, it can be shown that the estimate:
\begin{equation}
    \frac{d\mathcal{S}[\bm{u}]}{dt} \leq 0, 
    \label{eq:dsdtcontinuous}
\end{equation}
follows from \eqref{eq:ineqentconslaw} on periodic $\Omega$. This inequality will be the main interest of this paper.

\subsection{Entropy stable spatial discretization (FOM)}
We will discretize the conservation law \eqref{eq:conslaw} with a  finite volume method (FVM) based on flux-differencing \cite{leveque}. To introduce the general and frequently recurring structure of our entropy stability proof we will provide some detail on the full order model (FOM) discretization. The FVM discretization will be constructed such that discrete analogues to \eqref{eq:dsdtcontinuous} hold. Other discretization methods that similarly mimic \eqref{eq:dsdtcontinuous} are also possible, for example the split-form discontinuous Galerkin (DG) methods described in \cite{gassnersplit}, the summation-by-parts schemes in \cite{ranochacomparison} and the higher-order methods of \cite{fjordholmarbitrarily, fisherhigh}. We choose the FVM to keep the exposition simple, but note that our entropy-stable ROM framework should work with other entropy-stable FOM discretizations as well.

The scheme is formulated as:
\begin{equation}
    \Delta x_i \frac{d \bm{u}_i}{dt} + \bm{f}_{i+1/2} - \bm{f}_{i-1/2} = 0, \quad i \in \{0,...,N-1\},
    \label{eq:discconslaw}
\end{equation}
on a grid of $N$ grid cells so that $i \in \{0,...,N-1\}$. Here, $\bm{u}_i : [0,\infty) \rightarrow \mathbb{R}^n$ is the numerical solution vector in the $i$-th grid cell, $\Delta x_i := x_{i+1/2} - x_{i-1/2}$ is the cell size of the $i$-th cell with $x_{i\pm 1/2}$ denoting respectively the $x$ values of the right ($+$) and left ($-$) cell boundary and $\bm{f}_{i+1/2} := \bm{f}_h(\bm{u}_{i+1}, \bm{u}_{i})$ is the numerical flux on the right cell boundary of the $i$-th cell with $\bm{f}_h : \mathbb{R}^{n} \times \mathbb{R}^n \rightarrow \mathbb{R}^n$ being a two-point numerical flux function \cite{tadmornumerical} approximating the flux function $\bm{f}$ on a cell boundary based on two neighbouring numerical solution values, similarly $\bm{f}_{i-1/2}$ approximates $\bm{f}$ on the left boundary. We also define the total number of unknowns as $N_h := n\cdot N$. Periodic boundary conditions are enforced by setting $\boldsymbol{u}_{N} := \boldsymbol{u}_0$ and $\boldsymbol{u}_{-1} := \boldsymbol{u}_{N-1}.$ In the schemes we are considering the numerical flux function is constructed from entropy-conservative flux functions \cite{tadmorentropy, tadmornumerical, tadmorperfect}. These are flux functions that assure discrete analogues of \eqref{eq:entconslaw} and \eqref{eq:dsdtcontinuous} are satisfied with equality. This makes them suitable starting points from which to construct flux functions that have appropriate entropy-dissipative properties. We follow Tadmor's framework \cite{tadmornumerical} of entropy-conservative fluxes, which are defined as follows:
\begin{definition}[Entropy-conservative numerical flux]
    An entropy-conservative two-point numerical flux $\bm{f}^*_h : \mathbb{R}^n \times \mathbb{R}^n \rightarrow \mathbb{R}^n$ is a numerical two-point flux satisfying:
    \begin{enumerate}
        \item consistency: $\bm{f}^*_h(\boldsymbol{u},\boldsymbol{u}) = \bm{f}(\boldsymbol{u})$;
        \item symmetry: $\bm{f}_h^*(\boldsymbol{u}_l,\boldsymbol{u}_r) = \bm{f}_h^*(\boldsymbol{u}_r,\boldsymbol{u}_l)$;
        \item entropy conservation: $(\bm{\eta}(\boldsymbol{u}_l) - \bm{\eta}(\boldsymbol{u}_r))^T\bm{f}_h^*(\boldsymbol{u}_l,\boldsymbol{u}_r) = \psi(\boldsymbol{u}_l) - \psi(\boldsymbol{u}_r)$.
    \end{enumerate}
\end{definition}

The entropy-dissipative fluxes $\bm{f}_h : \mathbb{R}^n \times \mathbb{R}^n \times \mathbb{R}^{N_h} \rightarrow \mathbb{R}^n$ are now constructed from entropy-conservative fluxes $\bm{f}_h^*$ by adding (possibly solution dependent) entropy dissipation operators like:
\begin{equation}
    \bm{f}_{i+1/2} := \bm{f}_h(\boldsymbol{u}_{i+1}, \boldsymbol{u}_i, \boldsymbol{u}_h) = \bm{f}_h^*(\boldsymbol{u}_{i+1}, \boldsymbol{u}_i) - \bm{D}_{i+1/2}(\boldsymbol{u}_h)\Delta\bm{\eta}_{i+1/2}
    \label{eq:flux}
\end{equation}
with $\bm{D}_{i+1/2} : \mathbb{R}^{N_h} \rightarrow \mathbb{S}^n_+$ so that $\bm{D}_{i+1/2}(\boldsymbol{u}_h)$ is symmetric positive semi-definite (SPSD) for any $\boldsymbol{u}_h(t) \in \mathbb{R}^{N_h}$ ($\mathbb{S}^n_+$ and $\mathbb{S}^n_{++}$ are the convex sets of symmetric positive definite and symmetric positive semi-definite $n\times n$ matrices, respectively). Here, $\boldsymbol{u}_h : [0,\infty) \rightarrow \mathbb{R}^{N_h}$ is the numerical solution vector on the whole grid to be defined in what follows (this is required for higher order reconstructions like in e.g.\ \cite{fjordholmarbitrarily}). Additionally, we have defined $\Delta \bm{\eta}_{i+1/2} := \bm{\eta}(\boldsymbol{u}_{i+1}) - \bm{\eta}(\boldsymbol{u}_i)$. We will refrain from denoting explicitly the dependence on $\boldsymbol{u}_h$ in the third argument of $\bm{f}_h$ and simply write $\bm{f}_h(\boldsymbol{u}_{i+1}, \boldsymbol{u}_i)$ for \eqref{eq:flux}. 

For the purpose of model reduction in section \ref{sec:rom} we rewrite discretization \eqref{eq:discconslaw} with flux \eqref{eq:flux} in a matrix-vector formulation. We will introduce the following notations: volume-based quantities which live on cell centers and interface-based quantities which live on cell interfaces. The numerical solution vector is a volume-based quantity given by:
\begin{equation*}
    \bm{u}_h(t) := [u^1_0(t)...,u^k_i(t),u^k_{i+1}(t),...,u^n_i(t),...]^T \in \mathbb{R}^{N_h},
\end{equation*}
where $u^k_i(t) \in \mathbb{R}$ is the approximation of the $k$-th conserved variable (the variable conserved by the $k$-th equation in \eqref{eq:conslaw}) in the $i$-th cell evaluated at time $t$. The numerical flux vector is an interface-based quantity, overloading the notation for the numerical flux function, it is given by:
\begin{equation*}
    \bm{f}_h(\bm{u}_h) := [f^1_{1/2}, ..., f^1_{N-1/2},f^2_{1/2}, ... f^k_{i-1/2}, f^k_{i+1/2},...,f^n_{i+1/2},...,f^n_{N-1/2}]^T \in \mathbb{R}^{N_h},
\end{equation*}
where $f^k_{i+1/2}$ is the numerical flux of the $k$-th conservation equation evaluated at interface $i+1/2$ between cells $i+1$ and $i$. Periodic boundary conditions are built into the flux vector by evaluating $f^k_{N-1/2}(\bm{u}_0,\bm{u}_{N-1})$ for all $k\in\{1,...,n\}$.

To perform finite-difference operations as in \eqref{eq:discconslaw} for volume-based and interface-based quantities, respectively, the following matrices are defined:
\begin{equation*}
    \bar{\Delta}_v = \begin{bmatrix} 
        1 & 0 & 0 & - 1 \\
        -1 & 1 & 0 & 0 \\
        0 & \ddots & \ddots & 0 \\
        0 & 0 & -1 & 1 \\
    \end{bmatrix} \in \mathbb{R}^{N\times N}, \quad 
    \bar{\Delta}_i = \begin{bmatrix} 
        -1 & 1 & 0 & 0 \\
        0 & \ddots & \ddots & 0 \\
        0 & 0 & -1 & 1 \\
        1 & 0 & 0 & -1 \\
    \end{bmatrix} \in \mathbb{R}^{N\times N}.
\end{equation*}
We note the skew-adjointness relation
\begin{equation}
    \bar{\Delta}_v = -\bar{\Delta}_i^T, \label{eq:skew}
\end{equation}
and that both have zero row and column sum. These properties will be used in proving entropy stability of the scheme. For systems we define $\Delta_v := I \otimes \bar{\Delta}_v$ and $\Delta_i := I \otimes \bar{\Delta}_i$ with $I$ being the $n\times n$ identity matrix and $\otimes$ the Kronecker product. Clearly, $\Delta_v$ and $\Delta_i$ satisfy a similar skew-adjointness property \eqref{eq:skew}. We will also introduce the FVM mass matrices $\bar{\Omega}_h = \text{diag}(\Delta x_i)$ with $i = 0,1,...,N-1$ and $\Omega_h := I \otimes \bar{\Omega}_h$. With these operators we can write a compact form of the discretization \eqref{eq:discconslaw} as follows:
\begin{equation}
    \Omega_h\frac{d \bm{u}_h}{dt} + \Delta_v \bm{f}_h(\bm{u}_h) = 0.
    \label{eq:discconsmisc}
\end{equation}
To emphasize the role played by the dissipation operator $\bm{D}_{i+1/2}$ in obtaining entropy-stable spatial discretizations we will decompose $\Delta_v \bm{f}_h(\bm{u}_h)$ in an entropy-conserving part and an entropy-dissipating part, resulting in:
\begin{equation}
    \Omega_h\frac{d \bm{u}_h}{dt} + \Delta_v \bm{f}_h^*(\bm{u}_h) = \Delta_v \bm{D}_h(\bm{u}_h) \Delta_i \bm{\eta}_h,
    \label{eq:discretization}
\end{equation}
here, $\bm{f}_h^*(\bm{u}_h)$ is a vector of entropy conservative numerical fluxes, $\bm{D}_h(\bm{u}_h) \in \mathbb{S}_+^{N_h}$ is an SPSD matrix containing the terms associated to the dissipation operators $\bm{D}_{i+1/2}$ and $\bm{\eta}_h$ is a vector containing the grid values of the entropy variables ordered similarly as $\bm{u}_h$.

Having introduced an entropy-dissipative numerical flux, we evaluate the discrete analogue to the continuous total entropy functional \eqref{eq:entropycont} which should be suitably dissipated by the entropy stable discretization \eqref{eq:discretization} or conserved in the case of no dissipation. The discrete total entropy functional will be defined as:
\begin{equation}\label{eq:FOM_entropy}
    S_h[\bm{u}_h] := \bm{1}^T\bar{\Omega}_h \bm{s}_h,
\end{equation}
where $\bm{1}$ is a vector of ones and the local entropy is defined as:
\begin{equation*}
    \bm{s}_h(\bm{u}_h(t)) := [s(\bm{u}_0(t)), ..., s(\bm{u}_i(t)), ..., s(\bm{u}_{N-1}(t))]^T \in \mathbb{R}^N,
\end{equation*}
which is a volume-based quantity. The time evolution of $S_h$ is given by
\begin{equation*}
    \frac{dS_h[\bm{u}_h]}{dt} = \bm{1}^T\bar{\Omega}_h\frac{d\bm{s}_h}{dt} = \sum_i \Delta x_i \bm{\eta}(\bm{u}_i)^T \frac{d\bm{u}_i}{dt} = \bm{\eta}_h^T\Omega_h\frac{d\bm{u}_h}{dt}. 
\end{equation*}
We also define the entropy flux potential vector as:
\begin{equation*}
    \bm{\psi}_h(\bm{u}_h(t)) := [\psi(\bm{u}_0(t)), ..., \psi(\bm{u}_i(t)), ..., \psi(\bm{u}_{N-1}(t))]^T \in \mathbb{R}^N,
\end{equation*}
which is a volume-based quantity, like the local entropy vector. To analyse the entropy evolution we can substitute the spatial discretization \eqref{eq:discretization} in the previous expression to obtain:
\begin{align}
    \frac{dS_h[\bm{u}_h]}{dt} &= \bm{\eta}_h^T\Omega_h\frac{d\bm{u}_h}{dt} \nonumber \\
    &= -\bm{\eta}_h^T\Delta_v \bm{f}_h^*(\bm{u}_h) + \bm{\eta}_h^T\Delta_v\bm{D}_h(\bm{u}_h)\Delta_i\bm{\eta}_h \nonumber \\
    &= (\Delta_i \bm{\eta}_h)^T \bm{f}_h^*(\bm{u}_h) - \bm{\eta}_h^T \Delta_i^T \bm{D}_h(\bm{u}_h) \Delta_i \bm{\eta}_h \nonumber \\
    &= \bm{1}^T\bar{\Delta}_i \bm{\psi}_h - \bm{\eta}_h^T \Delta_i^T \bm{D}_h(\bm{u}_h) \Delta_i \bm{\eta}_h \nonumber \\
    &= 0 - \bm{\eta}_h^T \Delta_i^T \bm{D}_h(\bm{u}_h) \Delta_i \bm{\eta}_h \nonumber \\
    &\leq 0 \label{eq:dShdt},
\end{align}
where we used the skew-adjointness property \eqref{eq:skew}, the entropy conservation condition of the numerical flux, positive-definiteness of the dissipation operator $\bm{D}_h(\bm{u}_h)$ (and thus of $\Delta_i^T \bm{D}_h(\bm{u}_h) \Delta_i$) and the zero column sum of $\bar{\Delta}_i$. Clearly, in case no entropy dissipation is added in the numerical flux \eqref{eq:flux}, equation \eqref{eq:dShdt} reduces to
\begin{equation*}
    \frac{dS_h[\bm{u}_h]}{dt} = 0.
\end{equation*}
We note that the inequality \eqref{eq:dShdt} allows for formal $L^p$-stability statements \cite{svardweak, dafermoshyperbolic, fjordholmarbitrarily}.

\subsection{The entropy-stable linear ROM of Chan \cite{chanentropy}}\label{subsec:chan}
The main aim of this work is to propose reduced order models (ROMs) that are a nonlinear generalization of the entropy stable ROM of Chan \cite{chanentropy}. To highlight key conceptual differences between the ROM in \cite{chanentropy} and ours, and to introduce the ROM methodology, we will briefly discuss the elements of Chan's ROM leading to its entropy stability. Classical reduced order models including \cite{chanentropy} make the assumption that the evolution of $\bm{u}_h(t)$ can be accurately approximated with elements from a linear space $\mathcal{V} \subset \mathbb{R}^{N_h}$ where $\text{dim}(\mathcal{V}) := r \ll N_h$ so that $\mathcal{V}$ can be referred to as low-dimensional \cite{kleinenergy, rowleymodel, kalashnikovastable, arunajatesanstable, baronereduced, carlbergefficient, grimbergmesh, carlbergmodel}. Classically, the subspace $\mathcal{V}$ is constructed using a truncated proper orthogonal decomposition (POD) based on snapshot data collected in a matrix $X \in \mathbb{R}^{N_h \times n_s}$ with $n_s \in \mathbb{N}$ the number of snapshots \cite{sandersenonlinearly, kleinenergy}.
The construction of the ROM in \cite{chanentropy} starts by defining the approximation $\bm{u}_h \approx \bm{u}_r := \Phi \bm{a}$ with $\bm{a}\in\mathbb{R}^r$ being generalized coordinates in $\mathcal{V}$ relative to the basis $\Phi$. Here, we assume $\Phi$ is orthogonal in the $\Omega_h$-weighted inner product, i.e.\ $\phi_i^T\Omega_h \phi_j = \delta_{ij}$ with $\phi_i$ the $i$-th column of $\Phi$ and $\delta_{ij}$ the Kronecker delta function. Then, the approximation is substituted in \eqref{eq:discretization} introducing a semi-discrete residual
which is set orthogonal to $\mathcal{V}$ by solving the Galerkin projected system:
\begin{equation}
    \frac{d\bm{a}}{dt} + \Phi^T\Delta_v\bm{f}_r^*(\bm{a}) = \Phi^T\Delta_v \bm{D}_r(\bm{a})\Delta_i\bm{\eta}_r(\bm{a}),
    \label{eq:channoproj}
\end{equation}
with $\bm{f}_r^*(\bm{a}) := \bm{f}_h^*(\Phi \bm{a})$, $\bm{D}_r(\bm{a}) := \bm{D}_h(\Phi\bm{a})$ and $\bm{\eta}_r(\bm{a}) := \bm{\eta}_h(\Phi\bm{a})$. Equation \eqref{eq:channoproj} defines a (linear) POD-Galerkin ROM \cite{kleinenergy, rowleymodel, kalashnikovastable, arunajatesanstable, baronereduced}. 

Similar to the FOM case, see equation \eqref{eq:FOM_entropy}, we can evaluate the evolution of the ROM total entropy. The reduced total entropy functional is defined as:
\begin{equation}\label{eq:entropy_linear_ROM}
    S_r[\bm{a}] := S_h[\Phi \bm{a}] = \bm{1}^T\bar{\Omega}_h\bm{s}_r(\bm{a}),
\end{equation}
with $\bm{s}_r(\bm{a}) := \bm{s}_h(\Phi \bm{a})$. The evolution of the reduced total entropy is:
\begin{equation*}
    \frac{dS_r[\bm{a}]}{dt} = \bm{1}^T\bar{\Omega}_h\frac{d\bm{s}_r}{dt} = \sum_i \Delta x_i \bm{\eta}(\Phi_i \bm{a})^T \Phi_i \frac{d\bm{a}}{dt} = \bm{\eta}_r^T\Omega_h\Phi\frac{d\bm{a}}{dt},  
\end{equation*}
with $\Phi_i \in \mathbb{R}^{n \times r}$ the rows of $\Phi$ approximating values in cell $i$. The entropy evolution of \eqref{eq:channoproj} is:
\begin{align*}
    \frac{dS_r[\bm{a}]}{dt} &= \bm{\eta}_r^T\Omega_h\Phi\frac{d\bm{a}}{dt} \\
    &= - \bm{\eta}_r^T \Omega_h \Phi \Phi^T \Delta_v\bm{f}_r^*(\bm{a}) + \bm{\eta}_r^T \Omega_h \Phi \Phi^T\Delta_v \bm{D}_r( \bm{a})\Delta_i\bm{\eta}_r \\
    &= -\tilde{\bm{\eta}}_r^T \Delta_v\bm{f}_r^*(\bm{a}) + \tilde{\bm{\eta}}_r^T \Delta_v \bm{D}_r(\bm{a})\Delta_i\bm{\eta}_r,
\end{align*}
where, since $\Phi \Phi^T \Omega_h$ defines an $\Omega_h$-orthogonal projection operator \cite{axlerlinear}, $\tilde{\bm{\eta}}_r := \Phi \Phi^T \Omega_h \bm{\eta}_r$ are the so-called projected entropy variables. It is unclear whether this expression is bounded. To solve this, Chan \cite{chanentropy} proposes a technique used earlier in DG finite element literature \cite{parsanientropy, chandiscretely, chanontheentropy}. Namely, the discretization is not evaluated at $\Phi \bm{a}$ but at the entropy projected state:
\begin{equation*}
    \tilde{\bm{u}}_r := \bm{u}(\Phi \Phi^T \Omega_h \bm{\eta}_r) = \bm{u}(\tilde{\bm{\eta}}_r),
\end{equation*}
where we have defined $\bm{u} : \mathbb{R}^{N_h} \rightarrow \mathbb{R}^{N_h}, \bm{\eta}_h \mapsto \bm{u}_h$ for notational convenience. Recall that the mapping $\bm{u}$ is indeed available since $\bm{\eta}$ is injective. In this case we have:
\begin{align}
    \frac{dS_r[\bm{a}]}{dt} &= -\tilde{\bm{\eta}}_r^T \Delta_v\bm{f}_h^*(\bm{u}(\tilde{\bm{\eta}}_r)) + \tilde{\bm{\eta}}_r^T \Delta_v \bm{D}_h(\tilde{\bm{u}}_r)\Delta_i\tilde{\bm{\eta}}_r \nonumber \\
    &= (\Delta_i \tilde{\bm{\eta}}_r)^T \bm{f}_h^*(\bm{u}(\tilde{\bm{\eta}}_r)) - \tilde{\bm{\eta}}_r^T \Delta_i^T \bm{D}_h(\tilde{\bm{u}}_r)\Delta_i\tilde{\bm{\eta}}_r \nonumber \\
    &= \bm{1}^T\bar{\Delta}_i\tilde{\bm{\psi}}_r - \tilde{\bm{\eta}}_r^T \Delta_i^T \bm{D}_h(\tilde{\bm{u}}_r)\Delta_i\tilde{\bm{\eta}}_r \nonumber \\ 
    &= 0 - \tilde{\bm{\eta}}_r^T \Delta_i^T \bm{D}_h(\tilde{\bm{u}}_r)\Delta_i\tilde{\bm{\eta}}_r \nonumber\\
    &\leq 0, \label{eq:dSrdtlinear}
\end{align}
so that we re-obtain an entropy estimate that mimics the FOM estimate \eqref{eq:dShdt}. Here, $\tilde{\bm{\psi}}_r$ is the entropy flux potential vector evaluated at the entropy projected state. A key difference with the DG literature \cite{parsanientropy, chandiscretely, chanontheentropy} is that in the ROM case the basis $\Phi$ is only constructed to resolve solutions present in the snapshot matrix $X$, whereas in DG the trial basis is able to approximate a larger subspace of the relevant PDE function spaces. As a result, the DG trial basis can resolve the entropy variables well, but this may not be the case for the reduced basis $\Phi$.   
To address this issue, \cite{chanentropy} builds the basis $\Phi$ from a set of augmented snapshots given (with some abuse of notation) by:
\begin{equation*}
    \tilde{X} = [X, \bm{\eta}(X)],
\end{equation*}
so that the projection of the entropy variables on the basis $\Phi$ is close to the identity. As we will explain in section \ref{sec:tse}, for our proposed nonlinear manifold ROMs, such a construction is not sufficient, and a new tangent space enrichment technique will be proposed to ensure the accuracy of the entropy projection.

\section{An entropy-stable nonlinear manifold Galerkin ROM}\label{sec:rom}
The solution manifolds of many hyperbolic conservation laws \eqref{eq:conslaw} have slow Kolmogorov $n$-width decay \eqref{eq:knW}. Hence, approximations using linear subspaces as in \cite{chanentropy} may require very large reduced space dimensions $r$ before they become accurate. This comes at the cost of their efficiency. For this reason ROMs built on nonlinear spaces endowed (at least locally) with a manifold structure have become a topic of interest \cite{romorexplicable,romornonlinear,leemodel,sharmasymplectic,buchfinksymplectic,buchfinkmodel,barnettneural,barnettquadratic,tencerenabling,cocolahyper,kimfast,diazfast, geelenoperator, jainquadratic}. To address the shortcomings of the linear subspaces employed in \cite{chanentropy} we will generalize this method to nonlinear reduced spaces, while keeping the entropy-stability property. We will be interested specifically in nonlinear subsets of $\mathbb{R}^{N_h}$ endowed with some inner product, instead of any abstract space. Therefore we will not be very rigorous about our use of the word manifold, following predominantly the treatise of \cite{leemodel} and standard multivariable calculus interpretations. For a rigorous treatment we suggest consulting the recent preprint \cite{buchfinkmodel}. We will give a brief description of nonlinear manifold ROMs and then propose our generalization of \cite{chanentropy}.

\subsection{Manifold Galerkin model reduction}
In constructing ROMs on nonlinear manifolds we make the assumption that for any $t\in [0,T]$ there are points $\boldsymbol{u}_r(t)$ on a low-dimensional submanifold $\mathcal{M} \subset \mathbb{R}^{N_h}$ that accurately approximate $\boldsymbol{u}_h(t)$. Here, we denote $r := \text{dim}(\mathcal{M})$ and the low-dimensionality of $\mathcal{M}$ implies that $r \ll N_h$. We will refer to the submanifold $\mathcal{M}$ as the reduced manifold. Instead of the classical affine reduced space parameterization seen in the previous section we will use nonlinear manifold parameterizations given as:
\begin{equation}
    \boldsymbol{u}_h(t) \approx \boldsymbol{u}_r(t) := \bm{\varphi}(\boldsymbol{a}(t)) \in \mathcal{M},
    \label{eq:romappr}
\end{equation}
where:
\begin{equation}
    \bm{\varphi} : \mathbb{R}^r \rightarrow \mathbb{R}^{N_h},
    \label{eq:parameterization}
\end{equation}
is assumed to be a smooth nonlinear injective function - at least when restricted to some subset $\mathcal{A} \subseteq \mathbb{R}^r$ of interest where the ROM will be well-defined. This means that $\bm{\varphi}(\mathbb{R}^r) = \mathcal{M}$ with a Jacobian $\bm{J} : \mathbb{R}^r \rightarrow \mathbb{R}^{N_h \times r}, \hbox{ } \bm{a} \mapsto \frac{\partial \bm{\varphi}}{\partial \bm{a}}(\bm{a})$ of full rank for any $\boldsymbol{a} \in \mathcal{A} \subseteq \mathbb{R}^r$, where $\boldsymbol{a} : [0,T] \rightarrow \mathbb{R}^r$ are generalized coordinates on the manifold $\mathcal{M}$. The function $\bm{\varphi}$ may be obtained in many ways: some examples are quadratic approximations \cite{barnettquadratic, sharmasymplectic, buchfinkapproximation, geelenoperator, jainquadratic} or neural networks \cite{romorexplicable, romornonlinear, leemodel, kimfast, tencerenabling,diazfast,barnettneural,cocolahyper}. We will propose a new method based on rational polynomials in section \ref{sec:rational}. In this section, we develop an entropy-stable ROM which is agnostic of the choice for $\bm{\varphi}$.

To construct a ROM we substitute \eqref{eq:romappr} into the FOM discretization \eqref{eq:discretization} so that after applying the chain rule we find the residual:
\begin{equation*}
    \bm{r}\left(\frac{d\bm{a}}{dt},\bm{a}\right) := \bm{J} \frac{d\bm{a}}{dt} + \Omega_h^{-1}\Delta_v\bm{f}_r^*(\bm{a}) - \Omega_h^{-1} \Delta_v \bm{D}_r(\bm{a}) \Delta_i \bm{\eta}_r(\bm{a}),
\end{equation*}
where we define $\bm{f}_r^*(\bm{a}) := \bm{f}_h^*(\bm{\varphi}(\bm{a}))$, $\bm{D}_r(\bm{a}) := \bm{D}_h(\bm{\varphi}(\bm{a}))$ and $\bm{\eta}_r(\bm{a}) := \bm{\eta}_h(\bm{\varphi}(\bm{a}))$. The ROM is defined by minimizing this residual in the $\Omega_h$-norm for $\frac{d\bm{a}}{dt}$ given some $\bm{a}$, this results in the ROM:
\begin{equation*}
    (\bm{J}^T\Omega_h\bm{J}) \frac{d\bm{a}}{dt} + \bm{J}^T \Delta_v \bm{f}_r^*(\bm{a}) = \bm{J}^T\Delta_v\bm{D}_r(\bm{a}) \Delta_i \bm{\eta}_r(\bm{a}),
\end{equation*}
which is indeed well-defined for $\bm{a}\in\mathcal{A}\subseteq\mathbb{R}^r$ since the Jacobian $\bm{J}$ is assumed to be of full-rank on the subset $\mathcal{A}$, making the mass matrix $(\bm{J}^T\Omega_h\bm{J})$ invertible \cite{axlerlinear}. In this nonlinear case the ROM is given by the coefficients of an orthogonal projection of the FOM on the tangent space of $\mathcal{M}$ defined by $T_{u_r} \mathcal{M} := \text{span}(\bm{J}(\bm{a}))$ with $\bm{a}$ such that $\bm{u}_r = \bm{\varphi}(\bm{a})$. This orthogonal projection is carried out using the $\Omega_h$-weighted Moore-Penrose pseudoinverse $\bm{J}^\dagger := (\bm{J}^T\Omega_h\bm{J})^{-1}\bm{J}^T\Omega_h$. Constructing a ROM by projecting the FOM on the tangent space instead of the reduced manifold itself will result in key differences in our approach compared to the linear case outlined in \cite{chanentropy}: in contrast to the linear case where $\bm{J} = \Phi$, $T_{u_r} \mathcal{M}$ and $\mathcal{M}$ are no longer the same space. We introduce $\bm{J}^+ = \bm{J}^\dagger \Omega_h^{-1} = (\bm{J}^T\Omega_h\bm{J})^{-1}\bm{J}^T$ and write the ROM in compact form:
\begin{equation}
    \frac{d\bm{a}}{dt} + \bm{J}^+ \Delta_v \bm{f}_r^*(\bm{a}) = \bm{J}^+\Delta_v\bm{D}_r(\bm{a}) \Delta_i \bm{\eta}_r(\bm{a}).
    \label{eq:galrom}
\end{equation}
\begin{rmk}\label{rmk:nondim}
    The choice of inner-product spaces for ROMs of hyperbolic systems has recently come into question \cite{parishimpact}. Indeed when $n > 1$ the norm $||\bm{J}\bm{\alpha}||_{\Omega_h}$ for $\bm{\alpha} \in \mathbb{R}^r$ is dimensionally inconsistent in general. It is shown in \cite{parishimpact, kalashnikovastable} that dimensionally consistent inner products that are more appropriate in some sense can improve robustness of the ROMs. For our approach however it will be important that the ROM is calculated with the same inner product as used to calculate $\frac{dS_h}{dt}$. Therefore, we will only deal with nondimensionalized conservation laws. Alternatively, our results can also be applied at an equation-by-equation basis at the cost of potentially introducing a larger number of generalized coordinates.
\end{rmk}

\begin{rmk}
    A popular approach to construct nonlinear manifold ROMs is the least squares Petrov-Galerkin (LSPG) method \cite{romornonlinear, leemodel}. Using this method a fully discrete residual is minimized. We have chosen not to use this method because we want to use the structure of our entropy stable FOM discretization in constructing entropy stable ROMs. The fully discrete residual minimization approach of LSPG makes it more difficult to apply this structure. 
\end{rmk}

\subsection{An entropy stable nonlinear manifold Galerkin ROM}\label{eq:ES_nonlinear_ROM}
The reduced total entropy functional of the nonlinear manifold ROM is now defined similarly to the linear case (equation \eqref{eq:entropy_linear_ROM}) as:
\begin{equation*}
    S_r[\bm{a}] := S_h[\bm{\varphi}(\bm{a})] = \bm{1}^T \bar{\Omega}_h \bm{s}_r(\bm{a}),
\end{equation*}
with $\bm{s}_r(\bm{a}) := \bm{s}_h(\bm{\varphi}(\bm{a}))$. The reduced total entropy evolution is given by:
\begin{equation}
    \frac{dS_r[\bm{a}]}{dt} = \bm{1}^T \bar{\Omega}_h \frac{d\bm{s}_r}{dt} = \sum_i \Delta x_i \bm{\eta}(\bm{\varphi}_i(\bm{a}))^T\bm{J}_i \frac{d\bm{a}}{dt} = \bm{\eta}_r^T\Omega_h \bm{J}\frac{d\bm{a}}{dt},
    \label{eq:entropyprod}
\end{equation}
where $\bm{\varphi}_i : \mathbb{R}^r \rightarrow \mathbb{R}^n$ is the ROM approximation of the conserved variables in the $i$-th cell and $\bm{J}_i \in \mathbb{R}^{n\times r}$ is the Jacobian matrix of $\bm{\varphi}_i$ evaluated at $\bm{a}$. Using \eqref{eq:galrom}, the entropy evolution of the nonlinear manifold Galerkin ROM is:
\begin{align}
    \frac{dS_r[\bm{a}]}{dt} &= \bm{\eta}_r^T \Omega_h \bm{J}\frac{d\bm{a}}{dt} \nonumber \\
    &= -\bm{\eta}_r^T \Omega_h \bm{J}\bm{J}^+ \Delta_v \bm{f}_r^*(\bm{a}) + \bm{\eta}_r^T \Omega_h \bm{J}\bm{J}^+\Delta_v\bm{D}_r(\bm{a}) \Delta_i \bm{\eta}_r(\bm{a}) \nonumber \\
    &= -\tilde{\bm{\eta}}_r^T \Delta_v \bm{f}_r^*(\bm{a}) + \tilde{\bm{\eta}}_r^T \Delta_v\bm{D}_r(\bm{a}) \Delta_i \bm{\eta}_r(\bm{a}), \label{eq:rentevo}
\end{align}
where the projected entropy variables are defined as $\tilde{\bm{\eta}}_r = (\Omega_h\bm{J}\bm{J}^+)^T \bm{\eta}_r = \bm{J} \bm{J}^\dagger \bm{\eta}_r$, where $\bm{J} \bm{J}^\dagger$ is an $\Omega_h$-orthogonal projection on $T_{u_r}\mathcal{M}$. It follows from \eqref{eq:rentevo} that the reduced total entropy evolution satisfies an equation that is quite similar to the total entropy evolution of the FOM. However, instead of the actual entropy variables $\bm{\eta}_r$ evaluated at the point $\bm{u}_r$ on the reduced manifold $\mathcal{M}$, the inner product is taken with the projected entropy variables $\tilde{\bm{\eta}}_r$. We would like to use the entropy conservation condition at this point to show that the inner product in \eqref{eq:rentevo} is zero, but this does not hold because $\bm{\eta}(\bm{u}_r)\neq \tilde{\bm{\eta}}_r$ in general. To solve this, we can instead use the invertibility of the entropy variables to find for what state $\tilde{\bm{u}}_r \in \mathbb{R}^{N_h}$ we do have $\bm{\eta}(\tilde{\bm{u}}_r) = \tilde{\bm{\eta}}_r$. If we evaluate our flux at this state instead we can invoke the entropy conservation condition of the numerical flux to complete the proof of entropy conservation (or stability). This is exactly what is done in the linear setting in \cite{chanentropy} and leads to our main novelty.


We now present the main novelty of our work. 
We introduce a novel nonlinear manifold generalization of the linear entropy projection of \cite{chanentropy}. It is given by:
\begin{equation}
    \tilde{\bm{u}}_r = \bm{u}(\bm{J}\bm{J}^\dagger\bm{\eta}_r) = \bm{u}(\tilde{\bm{\eta}}_r),
    \label{eq:entproj}
\end{equation}
where the entropy variables $\bm{\eta}_r$ evaluated at the ROM state $\bm{u}_r$ are projected onto the tangent space $T_{u_r}\mathcal{M}$ instead of the reduced space itself. We note that projecting the entropy variables on the tangent space $T_{u_r}\mathcal{M}$ is a rather natural operation as the entropy variables can be interpreted as the gradient vector field of the entropy functional $S_h$ and thus as tangent vectors of $\mathbb{R}^{N_h}$. The vector $\tilde{\bm{u}}_r$ is the entropy projected state of the ROM. Carrying out this entropy projection on the tangent space is necessary as the projected entropy variables $\tilde{\bm{\eta}}_r$ appearing in the reduced total entropy evolution equation also are projected on the tangent space. 

A potential issue of our proposed form \eqref{eq:entproj} is that the difference between the projected entropy variables $\tilde{\bm{\eta}}_r$ and the original entropy variables $\bm{\eta}_r$ can be very large. Namely, in similar fashion to the naively constructed linear spaces spanned by $\Phi$ described in \autoref{subsec:chan} and \cite{chanentropy}, for arbitrary $\bm{\varphi}$ the entropy variables $\bm{\eta}_r$ may not be well-resolved by the columns of $\bm{J}(\bm{a})$ and thus by the $\Omega_h$-orthogonal projection in \eqref{eq:entproj}. This will very likely cause problems with accuracy of the ROM and the mapping $\bm{u}$ might not even be well-defined at $\tilde{\bm{\eta}}_r$. We will propose a novel method to assure the difference between $\bm{\eta}_r$ and $\tilde{\bm{\eta}}_r$ remains small in the following section. 

\begin{figure}
    \centering
    \includegraphics[width = 0.8\textwidth]{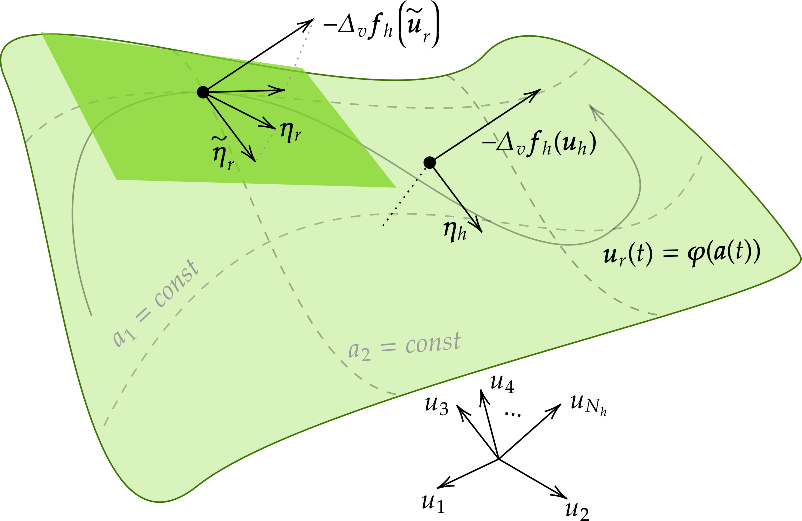}
    \caption{A visualization of the ROM construction. The entropy variables $\bm{\eta}_r$ are projected on the tangent space $T_{u_r}\mathcal{M}$ to obtain $\tilde{\bm{\eta}}_r$. A new state $\bm{u}_h$ (not necessarily on the reduced manifold) is found with entropy variables $\bm{\eta}_h$ such that $\tilde{\bm{\eta}}_r = \bm{\eta}_h$. We set $\tilde{\bm{u}}_r = \bm{u}_h$ and project $-\Delta_v \bm{f}_h(\tilde{\bm{u}}_r)$ orthogonally on the tangent space to complete the ROM.}
    \label{fig:manif}
\end{figure}

First, we continue constructing an entropy stable nonlinear manifold ROM and perform a Galerkin projection of the FOM at the entropy projected state $\tilde{\bm{u}}_r$. Doing so we obtain the following ROM:
\begin{equation}
    \frac{d\bm{a}}{dt} + \bm{J}^+\Delta_v\bm{f}_h^*(\tilde{\bm{u}}_r) = \bm{J}^+\Delta_v D_h(\tilde{\bm{u}}_r) \Delta_i \tilde{\bm{\eta}}_r,
    \label{eq:esROM}
\end{equation}
where we used $\tilde{\bm{\eta}}_r = \bm{\eta}_h(\bm{u}(\tilde{\bm{\eta}}_r))$. Here, the Jacobian matrix is still evaluated at the reduced coordinate $\bm{a} \in \mathbb{R}^r$ such that $\bm{u}_r = \bm{\varphi}(\bm{a})$ i.e. the non-projected state. It can be seen that the reduced total entropy evolution is bounded for this ROM since:
\begin{align}
    \frac{dS_r[\bm{a}]}{dt} &= -\tilde{\bm{\eta}}_r^T\Delta_v \bm{f}_h^*(\bm{u}(\tilde{\bm{\eta}}_r)) + \tilde{\bm{\eta}}_r^T\Delta_v\bm{D}_h(\tilde{\bm{u}}_r)\Delta_i\tilde{\bm{\eta}}_r \nonumber \\
    &= (\Delta_i \tilde{\bm{\eta}}_r)^T\bm{f}_h^*(\bm{u}(\tilde{\bm{\eta}}_r)) - \tilde{\bm{\eta}}_r^T\Delta_i^T\bm{D}_h(\tilde{\bm{u}}_r)\Delta_i\tilde{\bm{\eta}}_r \nonumber \\
    &= \bm{1}^T\bar{\Delta}_i\tilde{\bm{\psi}}_r - \tilde{\bm{\eta}}_r^T\Delta_i^T\bm{D}_h(\tilde{\bm{u}}_r)\Delta_i\tilde{\bm{\eta}}_r \nonumber \\
    &= 0 - \tilde{\bm{\eta}}_r^T\Delta_i^T\bm{D}_h(\tilde{\bm{u}}_r)\Delta_i\tilde{\bm{\eta}}_r \nonumber \\
    &\leq 0, \label{eq:romesstability}
\end{align}
where $\tilde{\bm{\psi}}_r = \bm{\psi}_h(\tilde{\bm{u}}_r)$ is the entropy flux potential of the entropy projected state and we were allowed to invoke the entropy conservation condition of the numerical fluxes. Clearly, we have:
\begin{equation*}
    \frac{dS_r[\bm{a}]}{dt} = 0,
\end{equation*}
when no entropy dissipation is present. Note that this approach exactly recovers the linear approach of \cite{chanentropy} when $\bm{J} = \Phi$, making it a proper generalization. Thus, by changing the state at which the FOM is evaluated and projected to the entropy projected state $\tilde{\bm{u}}_r$, correct total entropy evolution estimates can be recovered. We added a visualization of the ROM construction with an entropy projection in \autoref{fig:manif}.

\begin{rmk}
    To make relation \eqref{eq:romesstability} hold in a fully-discrete setting for Galerkin ROMs, entropy stable time integration is necessary. This is not trivial as most methods to satisfy entropy inequalities exactly during time integration require convexity of the entropy for existence results \cite{ranocharelaxation, tadmorentropy, ketchesonrelaxation}. Although this is the case for $S_h$, $S_r$ is not necessarily convex in the generalized coordinates $\bm{a}$ unless $\bm{\varphi}(\bm{a})$ is affine i.e.\ $\bm{\varphi}(\bm{a}) = \Phi \bm{a} + \bm{u}_0$ for some constant $\bm{u}_0 \in \mathbb{R}^{N_h}$ and $\Phi \in \mathbb{R}^{N_h\times r}$. 
    In this work, we focus on semi-discrete entropy stability and leave fully discrete entropy-stable ROMs as a suggestion for future work. In numerical experiments we use sufficiently small time steps to make entropy errors coming from the time integration negligible, for details see \autoref{sec:shallow}.
\end{rmk}

\begin{rmk}
    The entropy conservative hyper-reduction method proposed in \cite{chanentropy} does not generalize to nonlinear spaces as it relies on precomputation of compositions of linear operators. In the nonlinear case precomputation is not possible due to the changing tangent space. An entropy conservative hyper-reduction method suitable for nonlinear reduced spaces is also a suggestion for future work.
\end{rmk}

\subsection{Tangent space enrichment}\label{sec:tse}
To arrive at the correct entropy estimate \eqref{eq:romesstability} we carried out the entropy projection \eqref{eq:entproj}. Though the evolution of the entropy then satisfies a correct estimate, it is not clear whether the ROM solution itself remains accurate. Particularly, the difference between the entropy projected state \eqref{eq:entproj} and the original state \eqref{eq:romappr} can be very large. To see this we consider the entropy projection error:
\begin{equation*}
    \varepsilon_s := ||\bm{u}_r - \tilde{\bm{u}}_r||_{\Omega_h}.
\end{equation*}
Assuming the mapping $\bm{u}$ from entropy variables to conservative variables is sufficiently smooth, using the mean-value theorem we can bound this term as follows:
\begin{align*}
    ||\bm{u}_r - \tilde{\bm{u}}_r||_{\Omega_h} &= ||\bm{u}_r - \bm{u}(\bm{J}\bm{J}^\dagger\bm{\eta}_r)||_{\Omega_h} \\
    &= ||\bm{u}_r - \left[\bm{u}_r - \frac{\partial \bm{u}}{\partial \bm{\eta}}(\bm{\theta}) (I-\bm{J}\bm{J}^\dagger)\bm{\eta}_r\right]||_{\Omega_h} \quad \theta_i \in [(\bm{u}_r)_i, (\tilde{\bm{u}}_r)_i] \hbox{ } \forall i \\
    &\leq ||\frac{\partial \bm{u}}{\partial \bm{\eta}}(\bm{\theta})||_{\Omega_h} ||(I-\bm{J}\bm{J}^\dagger)\bm{\eta}_r||_{\Omega_h} \\
    &= ||\left(\frac{\partial^2 \bm{s}}{\partial \bm{u}^2}\right)^{-1}(\bm{\theta})||_{\Omega_h} ||(I-\bm{J}\bm{J}^\dagger)\bm{\eta}_r||_{\Omega_h},
\end{align*}
where we used the definition of the entropy variables in the last line and used the induced operator $\Omega_h$-norm for the Hessian matrix of the entropy function. It can be seen that there are two contributions to this bound. There is one based on the model, specifically on the Hessian of the entropy, and one given by the projection error of $\bm{\eta}_r$ on the tangent space $T_{u_r}\mathcal{M}$. As long as the entropy is a convex function at the mean value $\bm{\theta}$, the contribution of the model-based term is bounded by a term involving the inverse of the smallest eigenvalue of the entropy Hessian. We have little influence over this term. Specifically, this term can be large when the entropy is close to being non-convex. We do have control over the projection error. The magnitude of this term is controlled by the choice of reduced space. For a general reduced space constructed to contain solution snapshots this term can be very large, since the columns of the Jacobian $\bm{J}$ can be close to orthogonal to $\bm{\eta}_r$ while the standard nonlinear manifold ROM \eqref{eq:galrom} works fine. In the linear case, Chan \cite{chanentropy} solved this problem by enriching the snapshot data to construct $\Phi$ with snapshots of the entropy variables. This lowered the projection error contribution to the bound on $\varepsilon_s$ since $\bm{J}=\Phi$ in this case. However, for the general nonlinear case, $\bm{J}$ is not the same as the reduced space itself and we can no longer construct our reduced space to also contain the entropy variables to keep the projection error low. Instead, we need a different approach and therefore we propose a novel method to which we refer as \textit{tangent space enrichment}.

The key idea of tangent space enrichment is to construct an $r+1$-dimensional manifold $\hat{\mathcal{M}} \subset \mathbb{R}^{N_h}$ from the original $r$-dimensional manifold $\mathcal{M}$ by a `lifting' operation. Consequently, we use this new manifold for the ROM instead. This lifting operation is defined so that the original manifold $\mathcal{M}$ is a subset of the lifted manifold $\hat{\mathcal{M}}$, i.e.\ $\mathcal{M} \subset \hat{\mathcal{M}}$. Most importantly however, for all points $\bm{u}_r \in \mathcal{M} \subset \hat{\mathcal{M}}$ \textit{the lifting operation is constructed such that} $\bm{\eta}_r \in T_{u_r}\hat{\mathcal{M}}$. This means that the entropy variable projection error is precisely zero at the points contained in the \textit{old} manifold when projecting $\bm{\eta}_r$ on the tangent space of the \textit{new} manifold. This assures that the entropy projection is accurate for the points $\bm{u}_r \in \mathcal{M} \subset \hat{\mathcal{M}}$ when using tangent space enrichment. 

We motivate this approach over a more straightforward generalization of Chan's snapshots enrichment \cite{chanentropy} method by the following. Nonlinear reduced spaces are often constructed iteratively by minimizing some loss function. A nonlinear version of Chan's enrichment method would require, for a given $\bm{a}$, fitting $\bm{\varphi}(\bm{a})$ to a snapshot $\bm{u}_h$ whilst the Jacobian $\bm{J}(\bm{a})$ has a low entropy projection error $\varepsilon_s$. The construction of $\mathcal{M}$ would therefore require including terms based on $\bm{J}$ in the loss function. This can be very expensive and generally will not exactly embed the entropy variables in the tangent space at the appropriate points. As will be discussed, our method requires no extra effort in constructing $\bm{\varphi}$ and contrary to the straightforward approach exactly enriches the tangent spaces with the correct entropy variables.


Our novel tangent space enrichment is defined by the following parameterization for $\hat{\mathcal{M}}$: 
\begin{equation}
    \hat{\bm{\varphi}}(\bm{a},\alpha) = \bm{\varphi}(\bm{a}) + \bm{\eta}(\bm{\varphi}(\bm{a})) \alpha,
    \label{eq:tseparam}
\end{equation}
here $\bm{\varphi} : \mathbb{R}^r \rightarrow \mathbb{R}^{N_h}$ is the parameterization of the original manifold $\mathcal{M}$, $\bm{a} \in \mathbb{R}^r$ are the $r$ reduced coordinates associated to the original parameterization $\bm{\varphi}$ and $\alpha \in \mathbb{R}$ is the $r+1$-th reduced coordinate associated to the lifting operation. As in remark \autoref{rmk:nondim}, here we see another reason for the importance of non-dimensionalization. Namely, on dimensional grounds the expression \eqref{eq:tseparam} does not make sense if $\bm{\varphi}$ and $\bm{\eta}$ are not suitably normalized. 

\begin{figure}
    \centering
    \includegraphics[width = 0.8\textwidth]{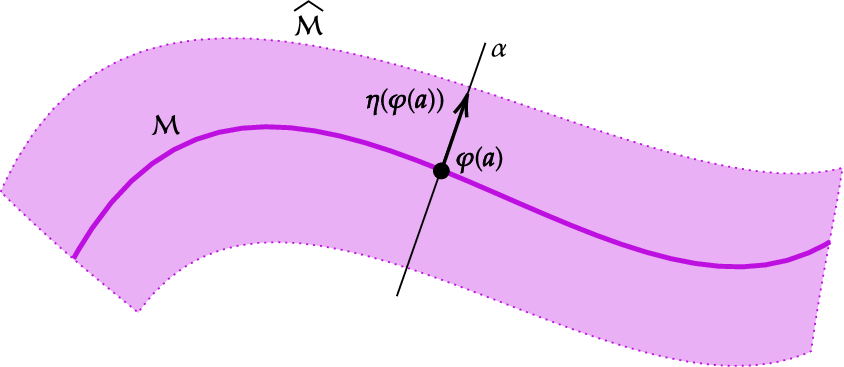}
    \caption{A visualization of tangent space enrichment. The dark purple curve is the original manifold $\mathcal{M}$, the light purple region is a section of the enriched manifold $\mathcal{M}^*$. $\mathcal{M}^*$ is constructed by lifting points from $\varphi(\bm{a})$ along lines in the direction of $\bm{\eta}(\varphi(\bm{a}))$.}
    \label{fig:tse}
\end{figure}

The new parameterization can be interpreted as follows. Given any point $\bm{\varphi}(\bm{a}) \in \mathcal{M}$ we generate new points $\hat{\bm{u}}_r \in \hat{\mathcal{M}}$ by lifting the new points from the point $\bm{\varphi}(\bm{a})$ in the direction of $\bm{\eta}(\bm{\varphi}(\bm{a}))$ by a distance $||\hat{\bm{\varphi}}(\bm{a},\alpha) - \bm{\varphi}(\bm{a})|| = ||\bm{\eta}(\bm{\varphi}(\bm{a}))|| \cdot |\alpha|$. Note that at $\alpha = 0$ we do not lift the point at all, resulting in a point at $\bm{\varphi}(\bm{a})$; in other words, the original manifold $\mathcal{M}$ is the set $\hat{\bm{\varphi}}(\mathbb{R}^r,\alpha=0)$. The lifting operation is visualized in \autoref{fig:tse}. The Jacobian matrix of the new parameterization \eqref{eq:tseparam}, whose columns span the tangent space $T_{\hat{u}_r}\hat{\mathcal{M}}$ of the lifted manifold $\hat{\mathcal{M}}$ at the point $\hat{\bm{u}}_r \in \hat{\mathcal{M}}$, is given by:
\begin{align}
    \hat{\bm{J}}(\bm{a},\alpha) &= \begin{bmatrix}
        \frac{\partial \hat{\bm{\varphi}}}{\partial \bm{a}} & \frac{\partial \hat{\bm{\varphi}}}{\partial \alpha}
    \end{bmatrix} \label{eq:tsejacobian}\\
    &= \begin{bmatrix}
        \left(I + \alpha  \hbox{ } \frac{\partial \bm{\eta}}{\partial \bm{u}}(\bm{\varphi}(\bm{a}))\right)\bm{J}(\bm{a}) & \bm{\eta}(\bm{\varphi}(\bm{a}))
    \end{bmatrix}, \nonumber
\end{align}
where $I \in \mathbb{R}^{N_h \times N_h}$ is the identity matrix and $\frac{\partial \bm{\eta}}{\partial \bm{u}}(\bm{\varphi}(\bm{a})) = \frac{\partial^2 \bm{s}_h}{\partial \bm{u}^2}(\bm{\varphi}(\bm{a})) \succeq 0$ is a sparse SPSD $2n-1$-diagonal matrix containing components of the Hessian of the local entropy value with respect to the solution on each diagonal. Note that the derivative with respect to the $r+1$-th tangent space enrichment coordinate $\alpha$ is exactly $\bm{\eta}(\bm{\varphi}(\bm{a}))$. Furthermore, on the old manifold $\mathcal{M}$ associated to $\alpha = 0$ the matrix $\frac{\partial \hat{\bm{\varphi}}}{\partial \bm{a}}$ is equal the original Jacobian $\bm{J}(\bm{a})$. At the points $\hat{\bm{\varphi}}(\bm{a},\alpha = 0) = \bm{\varphi}(\bm{a})$ we have thus exactly enriched the tangent space with the entropy variables $\bm{\eta}_r = \bm{\eta}(\bm{\varphi}(\bm{a}))$ at those points. This is the direct result of the lifting operation. This can be seen from the enriched parameterization \eqref{eq:tseparam}. Lifting a point from $\bm{\varphi}(\bm{a})$ by changing $\alpha$ while keeping $\bm{a}$ constant, moves a point in the direction tangent to $\bm{\eta}(\bm{\varphi}(\bm{a}))$. As a consequence $\bm{\eta}(\bm{\varphi}(\bm{a}))$ appears as a tangent vector in the enriched Jacobian \eqref{eq:tsejacobian}.

Given the tangent space enrichment, the ROM is constructed in a similar fashion as in the previous section, but using the enriched lifted manifold $\hat{\mathcal{M}}$. The $r+1$-th tangent space enrichment coordinate $\alpha$ is simply treated as an additional reduced coordinate. The ROM thus takes the form:
\begin{equation*}
    \frac{d}{dt}\begin{bmatrix}
        \bm{a} \\ \alpha
    \end{bmatrix} + \hat{\bm{J}}^+\Delta_v \bm{f}_h^*(\hat{\tilde{\bm{u}}}_r) = \hat{\bm{J}}^+ \Delta_v \bm{D}_h(\hat{\tilde{\bm{u}}}_r)\Delta_i \hat{\tilde{\bm{\eta}}}_r,
\end{equation*}
where:
\begin{equation*}
    \hat{\tilde{\bm{u}}}_r = \bm{u}(\hat{\bm{J}}\hat{\bm{J}}^\dagger\hat{\bm{\eta}}_r),
\end{equation*}
and $\hat{\tilde{\bm{\eta}}}_r = \hat{\bm{J}}\hat{\bm{J}}^\dagger \bm{\eta}(\hat{\bm{\varphi}}(\bm{a},\alpha))$. Since the proof of entropy stability for our ROM is independent of the manifold parameterization, this ROM is still entropy stable.

\section{Rational polynomial manifolds}\label{sec:rational}
\subsection{Background}
In this article we are interested in systems that can exhibit strong spatial gradients that are moving in time. Many existing data compression methods for reduced manifold construction are not well suited for these types of systems. Linear data compression methods like proper orthogonal decomposition (POD) and general reduced basis methods fail because the moving gradients imply that the data is often of very high rank. This indicates that the data is not well represented in low-dimensional linear subspaces, which is the fundamental assumption of linear approaches. Nonlinear data compression methods offer a potential solution to this problem by instead compressing the highly nonlinear data on nonlinear reduced manifolds. In model reduction different parameterizations have become popular, in particular the decoder part of autoencoder neural networks \cite{romorexplicable, romornonlinear, leemodel, kimfast, tencerenabling,diazfast,barnettneural,cocolahyper} and multivariate quadratic polynomials \cite{barnettquadratic, sharmasymplectic, buchfinkapproximation, geelenoperator, jainquadratic}. However, in the vicinity of strong gradients these nonlinear methods can suffer from oscillations \cite{geelenoperator}. These oscillations may be difficult or impossible to remove. This is a problem for Galerkin projection-based ROMs which can be sensitive to errors in the solution compression \cite{carlbergmodel}. To accurately assess the performance of our novel entropy stable manifold Galerkin ROM there is thus a need for nonlinear data compression methods that are more capable of dealing with large and moving spatial gradients, particularly without significant oscillations. Recently, neural networks with discontinuous activation functions have been proposed \cite{dellasantadiscontinuous}, however training these networks can be cumbersome. Furthermore, registration based approaches \cite{taddeiregistration, mojganilow} have been very effective, but have not yet been applied in the context of nonlinear manifold ROMs similar to ours. In this research we will propose a novel reduced manifold parameterization method based on rational polynomials, that is far less oscillatory around strong spatial gradients than the previously mentioned methods (neural networks and quadratic approaches), but is still equally interpretable as the recently proposed quadratic manifolds. 

\subsection{Pole-free rational quadratic manifolds}
We now give a description of rational polynomial manifolds. A rational polynomial manifold is the element-wise ratio of two polynomial manifolds:
\begin{equation}
    \bm{\varphi}(\bm{a}) = \frac{\sum_{i=1}^{p_{\text{num}}} \bm{H}^i : \bm{a}^{\otimes i} + \bm{u}_{\text{ref}}}{\sum_{i=1}^{p_{\text{den}}}\bm{G}^i : \bm{a}^{\otimes i} + \bm{1}},
    \label{eq:ratpolyman}
\end{equation}
here, $\bm{H}^i, \bm{G}^i \in \mathbb{R}^{N_h \times r \times ... \times r}$ are $(i+1)^{\text{th}}$-order tensors with the first axis of size $N_h$ and $i$ axes of length $r$, $\bm{a}^{\otimes i}$ is the $i$-fold outer product such that for example $(\bm{a}^{\otimes 3})_{ijk} = a_i a_j a_k$, $\bm{H}^i : \bm{a}^{\otimes i}$ denotes summation of the components of $\bm{a}^{\otimes i}$ and the components of slices along the first axis of $\bm{H}^i$, again as example $\left(\bm{H}^3 : \bm{a}^{\otimes 3}\right)_i = \sum_{j,k,l=0}^{r-1} \left(\bm{H}^3\right)_{ijkl} a_j a_k a_l$. Furthermore, we have $\bm{u}_{\text{ref}} \in \mathbb{R}^{N_h}$ and we consider division of two vectors element-wise. The constant vector in the denominator has been set to one without loss of generality. The expression \eqref{eq:ratpolyman} generalizes polynomial manifolds \cite{geelenoperator, barnettquadratic}, in the sense that those are recovered by setting $p_{\text{den}} = 0$. This shows that rational polynomial manifold encapsulate a larger class of functions than polynomial manifolds. 
By introducing a polynomial in the denominator and allowing it to rapidly and smoothly approach zero for small changes in $\bm{a}$ we can obtain very fast and smooth increases in the function value of $\varphi$ without oscillations. This gives us the opportunity to model steep gradients in the snapshot data which may be present in the form of advected shocks. Since higher-order tensors can become quite expensive to deal with, we restrict our attention to rational quadratic manifolds by setting $p_{\text{num}} = p_{\text{den}} = 2$.

To compromise between efficiency and expressiveness we will take $p_{\text{num}} = p_{\text{den}} = 2$. The $i$-th component of the vector-valued function output $\varphi_i(\bm{a})$ can then be written as:
\begin{equation}
    \varphi_i(\bm{a}) = \frac{\bm{a}^T\bm{H}_i^2 \bm{a} + \bm{H}_i^1\bm{a} + (\bm{u}_{\text{ref}})_i}{\bm{a}^T\bm{G}_i^2 \bm{a} + \bm{G}_i^1\bm{a} + 1},
    \label{eq:quadratslice}
\end{equation}
where $\bm{H}_i^2,\bm{G}_i^2 \in \mathbb{R}^{r\times r}$ are the $i$-th slices along the first axes of $\bm{H}^2$ and $\bm{G}^2$, respectively and $\bm{H}_i^1, \bm{G}_i^1 \in \mathbb{R}^{1\times r}$ are the $i$-th rows of $\bm{H}^1$ and $\bm{G}^1$, respectively. Since the matrices are only used in quadratic forms, we can, without loss of generality, assume $\bm{H}_i^2,\bm{G}_i^2$ to be symmetric. An obvious concern with expressions of this form is the occurrence of spurious poles, i.e.\ unwanted division by zero. We avoid this issue for the case $p_{\text{num}} = p_{\text{den}} = 2$ and all $\bm{a} \in \mathbb{R}^r$ by constraining the quadratic form in the denominator to be positive semi-definite and setting the linear term to zero:
\begin{equation*}
    \bm{G}_i^2 \succeq 0, \quad \bm{G}_i^1 = 0, \quad \forall i \in \{0,...,N_h\}.
\end{equation*}
It is easily seen that, in this case, the denominator is never less than 1. Consequently, spurious poles cannot occur for any real $\bm{a}$. Removing the linear term has not resulted in a significant loss in accuracy in our numerical experiments. The full $\bm{\varphi}$ is then given as:
\begin{equation}
    \bm{\varphi}(\bm{a}) = \frac{\bm{H}^2 : [\bm{a} \otimes \bm{a}] + \bm{H}^1 \bm{a} + \bm{u}_{\text{ref}}}{\bm{G} : [\bm{a} \otimes \bm{a}] + \bm{1}}, \quad \bm{G}_i \succeq 0 \hbox{ } \forall i,
    \label{eq:quadrat}
\end{equation}
since there is no linear term in the denominator we write $\bm{G}$ instead of $\bm{G}^2$. To construct manifold Galerkin ROMs we will require the Jacobian matrix of this expression - see equations \eqref{eq:galrom} and \eqref{eq:tsejacobian}. The Jacobian matrix is given by the following:
\begin{equation*}
    \frac{\partial \bm{\varphi}}{\partial \bm{a}} = \frac{2\bm{H}^2 \cdot \bm{a} + \bm{H}^1}{\left(\bm{G} : [\bm{a} \otimes \bm{a}] + \bm{1}\right) \otimes \bm{1}} - \left[\frac{\bm{H}^2 : [\bm{a} \otimes \bm{a}] + \bm{H}^1 \bm{a} + \bm{u}_{\text{ref}}}{\left(\bm{G} : [\bm{a} \otimes \bm{a}] + \bm{1}\right)^2} \otimes \bm{1}\right] \circ (2\bm{G} \cdot \bm{a}),
\end{equation*}
where the operation $2\bm{H}^2 \cdot \bm{a} \in \mathbb{R}^{N_h \times r}$ indicates slice-wise matrix multiplication, i.e.\ for the $i$-th row it holds that $(2\bm{H}^2 \cdot \bm{a})_i = 2\bm{H}_i^2 \bm{a}$, due to symmetry of $\bm{H}_i^2,\bm{G}_i$ the order of axes is not relevant, division of matrices is understood element-wise and $\circ$ is the Hadamard matrix multiplication operator. 

\subsection{Manifold construction}
We will determine the coefficient tensors in \eqref{eq:quadrat} from data. Like for quadratic manifolds \cite{barnettquadratic, geelenoperator}, it holds for rational quadratic manifolds that the coefficients in two different slices of the coefficient tensors are independent, as can be seen in \eqref{eq:quadratslice}. Consequently, the coefficients can be determined purely from the parametric and temporal behaviour of the data in the specific cell and solution variable associated to the $i$-th component $\varphi_i(\bm{a})$ of $\bm{\varphi}(\bm{a})$. The task of fitting a rational manifold thus reduces to fitting an expression \eqref{eq:quadratslice} to data for each cell and solution variable. 

In the spirit of quadratic manifolds we will compress the snapshot data defined by:
\begin{equation*}
    X = [\bm{u}_h(t^0), \bm{u}_h(t^1), ..., \bm{u}_h(t^{n_s - 1})] \in \mathbb{R}^{N_h \times n_s}, 
\end{equation*} 
where $n_s \in \mathbb{N}$ is the number of snapshots, as the coefficients of their projection on the first $r$ left singular vectors of $X$ given in $\Phi \in \mathbb{R}^{N_h \times r}$. Following this we define $A := (\Phi^T X)^T \in \mathbb{R}^{n_s \times r}$. We aim to find the coefficients $\bm{G}_i, \bm{H}_i^2, \bm{H}_i^1, (\bm{u}_{\text{ref}})_i$ such that for all $j = 0, ..., n_s-1$ we have $\varphi_i(\bm{a}_j) \approx X_{ij}$, where $\bm{a}_j$ is the $j$-th row of $A$ written as a column vector. Defining the $i$-th row of $X$ as $\bm{y}^i \in \mathbb{R}^{n_s}$, the optimization procedure for the coefficients $\bm{G}_i, \bm{H}_i^2, \bm{H}_i^1, (\bm{u}_{\text{ref}})_i$ is formulated as:
\begin{equation}
    \bm{G}_i, \bm{H}_i^2, \bm{H}_i^1, (\bm{u}_{\text{ref}})_i = \argmin_{\bm{G} \in \mathbb{S}^r_+ ,\bm{H}\in \mathbb{S}^r,\bm{h} \in \mathbb{R}^r, u \in \mathbb{R}} \left|\left|\bm{y}^i - \frac{(A\bm{H} \circ A)\bm{1} + A \bm{h} + u}{(A\bm{G} \circ A)\bm{1} + 1} \right|\right|^2,
    \label{eq:optrat}
\end{equation} 
though the sets $\mathbb{S}^r_+$ and $\mathbb{S}^r$ are convex subsets of $\mathbb{R}^{r\times r}$ the objective function is non-convex due to the division operation. 

A popular approach is to linearize this nonlinear optimization problem \cite{hokansonmultivariate, austinpractical} by multiplying with the denominator. This results in a convex semi-definite program which can be solved very efficiently with e.g.\ interior point methods. However, the optimum value of this linearized problem is generally not the same as the nonlinear problem \eqref{eq:optrat}. Iterative approaches that attempt to somehow refine the solution of the linearized problem and that can potentially be supplemented with our semi-definite constraint exist \cite{hokansonmultivariate, austinpractical, gustavsenimproving}. However, convergence to (local) minima of the nonlinear problem is generally not guaranteed, nor is convergence in general \cite{shinonconvergence, hokansonleast}. For this reason, we will stick to the fully nonlinear and expensive optimization problem \eqref{eq:optrat}. Nonetheless, we believe that the linearized approach holds significant promise and that it will be crucial for future work in order to scale the approach to larger meshes and datasets. Finally, we will implement the semi-definite constraint in this fully nonlinear setting by optimizing for the Cholesky decomposition of $\bm{G}_i = L_i L_i^T$ \cite{axlerlinear}:
\begin{equation}
    L_i, \bm{H}_i^2, \bm{H}_i^1, (\bm{u}_{\text{ref}})_i = \argmin_{L \in \mathbb{L}^r ,\bm{H}\in \mathbb{S}^r,\bm{h} \in \mathbb{R}^r, u \in \mathbb{R}} \left|\left|\bm{y}^i - \frac{(A\bm{H} \circ A)\bm{1} + A \bm{h} + u}{(AL \circ AL)\bm{1} + 1} \right|\right|^2,
    \label{eq:optrat2}
\end{equation} 
where $\mathbb{L}^r \subset \mathbb{R}^{r\times r}$ is the vector subspace of lower triangular $r\times r$ matrices. As initial guess we can either use $L_{i-1}, \bm{H}_{i-1}^2, \bm{H}_{i-1}^1, (\bm{u}_{\text{ref}})_{i-1}$ if available and corresponding to the same solution variable, or otherwise simply vectors or tensors consisting of ``ones". We will carry out the fitting procedure using the JAXFit package \cite{jaxfit} for GPU-accelerated nonlinear least squares solutions.

\section{Numerical experiments}\label{sec:experiments}
To show that our entropy stable manifold Galerkin ROMs satisfy appropriate semi-discrete entropy inequalities we will perform numerical experiments on a range of one-dimensional nonlinear conservation laws. We will carry out the experiments using the rational quadratic manifolds proposed in section \ref{sec:rational}. We will also compare the ability of the rational quadratic manifolds to compress convection dominated data to that of linear POD-based methods and quadratic manifolds \cite{barnettquadratic, jainquadratic, geelenoperator}. We do not compare against neural network based approaches \cite{romornonlinear, leemodel, barnettneural} since in our experience they struggle with approximating discontinuities, and require careful hyperparameter tuning to give reasonable results. The underlying FOM discretizations will be the existing TeCNO schemes of \cite{fjordholmarbitrarily}, so that we will only mention some aspects of the discretizations but for details we will refer to \cite{fjordholmarbitrarily}. We will start with the inviscid Burgers equation in \autoref{sec:burgers}, then we will treat the shallow water equations in \autoref{sec:shallow} and finally we will treat the compressible Euler equations with ideal thermodynamics in \autoref{sec:euler}. We test different aspects of the ROM with the different test cases, an overview of the different test purposes has been provided in \autoref{tab:purpose}.

\begin{table}[h!]
\centering
\begin{tabular}{ll}
\toprule
Experiment         & Purpose                   \\ 
\midrule
Inviscid Burgers   & Manifold accuracy         \\
Shallow water      & Entropy conservation properties   \\
Compressible Euler & Impact of entropy projection and tangent space enrichment\\ 
\bottomrule
\end{tabular}
\caption{Overview of numerical experiments.}
\label{tab:purpose}
\end{table}

The experiments have been implemented using the JAX library \cite{jaxgithub} in Python, which allows for automatic differentiation to compute Jacobian matrices and, where possible, accelerated computing using an Nvidia A2000 laptop GPU.

\subsection{Inviscid Burgers equation}\label{sec:burgers}
We will use this experiment mainly to highlight our proposed rational quadratic manifolds when compared to existing (`standard') Galerkin ROMs on different types of manifolds. We will already include the entropy stable ROM \eqref{eq:esROM} here, but the focus on the role of entropy stability will be in the next test cases. 

The inviscid Burgers equation is given by:
\begin{equation}
    \frac{\partial u}{\partial t} + \frac{\partial}{\partial x}\left(\frac{u^2}{2}\right) = 0,
    \label{eq:burgers}
\end{equation}
with conserved variable $u : \Omega \times [0,T] \rightarrow \mathbb{R}$. As continuous and discrete entropy we take \cite{leveque}:
\begin{equation*}
    \mathcal{S}[u] = \frac{1}{2}\int_{\Omega} u^2 dx,
\end{equation*}
and
\begin{equation*}
    S_h[\bm{u}_h] = \frac{1}{2}||\bm{u}_h||_{\Omega_h}^2.
\end{equation*}
The reduced entropy functional follows in a straightforward fashion from the discrete entropy functional. Using these specific entropies we have for the entropy variables:
\begin{equation*}
    \eta(u) = u.
\end{equation*}
Consequently, the manifold parameterization with TSE is given in the particularly simple form:
\begin{equation*}
    \bm{\varphi}^*(\bm{a},\alpha) = (1 + \alpha) \cdot \bm{\varphi}(\bm{a}).
\end{equation*}
An entropy conservative flux is given by \cite{fjordholmarbitrarily}:
\begin{equation*}
    f_{i+1/2} = \frac{u_{i+1}^2 + u_{i+1}u_i + u_i^2}{6},
\end{equation*}
and we use a local Lax-Friedrichs-type of entropy dissipation operator \cite{leveque}:
\begin{equation*}
    D_{i+1/2}(\bm{u}_h) = \max(|u_{i+1}|,|u_i|).
\end{equation*}
We will discretize \eqref{eq:burgers} on a domain $\Omega = \mathbb{T}([0,L])$ of length $L = 1$ using a numerical grid consisting of 300 cells. Discretization in time will be done using the classical Runge-Kutta 4 (RK4) method \cite{griffithsnumerical} with a time step of size $\Delta t = 0.001$. We will integrate in time until $T = 1$. Solution snapshots are captured after every 5 timesteps resulting in $n_s = 201$. We perform two tests: in the first, we compress the data to a reduced dimension of $r = 15$ for all manifolds, and in the second we compress the data to a reduced dimension such that the reconstruction errors (to be defined later) are of the same order as the rational quadratic manifold with $r = 15$. The initial condition is a simple offset sine wave:
\begin{equation*}
    u_0(x) = \sin(2\pi x) + 1.
\end{equation*}

\begin{figure}
    \centering
    \includegraphics[width = 0.5\textwidth]{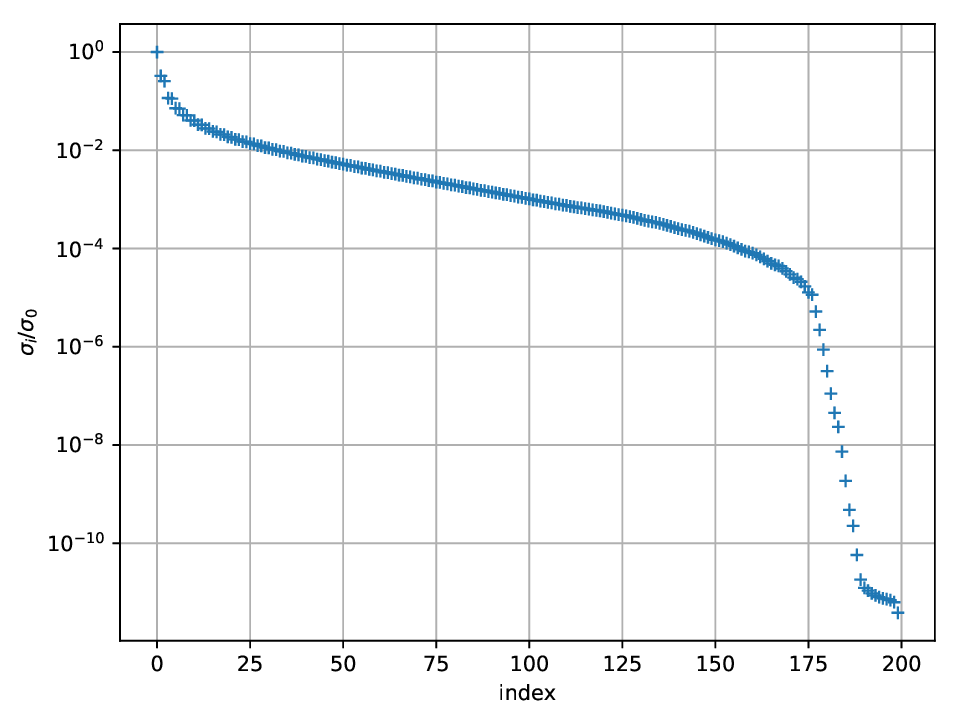}
    \caption{Normalized singular values of the inviscid Burgers equation data.}
    \label{fig:burgerssingvals}
\end{figure}

The KnW decay \eqref{eq:knW} for this system is very slow as is evident from the normalized singular values depicted in \autoref{fig:burgerssingvals}. Defining the relative information content (RIC) as in \cite{lassilamodel} we have $\text{RIC}\approx 99.5\%$ for $r = 15$. 

\begin{figure}
    \centering
    \includegraphics[width = 0.8\textwidth]{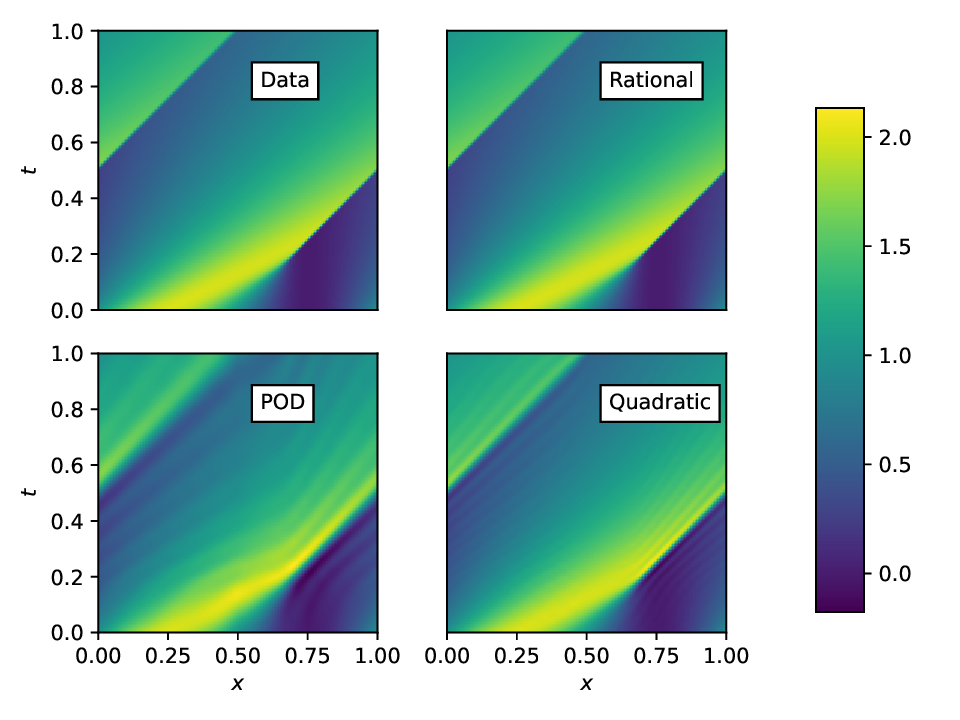}
    \caption{Original data and reconstructions for $r=15$ in space-time plots.}
    \label{fig:burgersreconstr}
\end{figure}

\begin{figure}
    \centering
    \includegraphics[width = \textwidth]{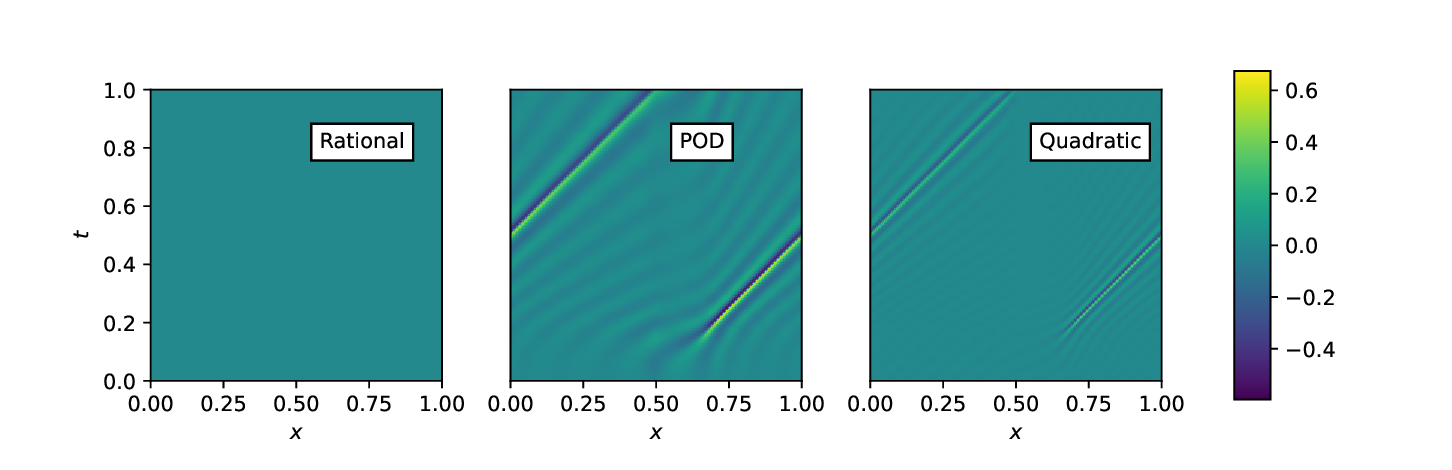}
    \caption{Space--time errors $\varepsilon_{xt}$ of the reconstructions.}
    \label{fig:burgersreconstrerror}
\end{figure}

We will first compare the reconstruction accuracy of our proposed rational quadratic manifold to existing quadratic manifold \cite{barnettquadratic} and POD linear \cite{sandersenonlinearly} manifold approaches for the solution data of the FOM with $r = 15$. To do this we will save the matrix of generalized coordinates $A = \Phi^TX \in \mathbb{R}^{r\times n_s}$ associated to the snapshots in $X$. Note that these coordinates form the reduced representation for all manifolds since all of the tested manifolds are constructed from the POD compression of the data. In turn we will try to reconstruct the snapshots in $X$ from their reduced representations in $A$. We will construct the quadratic manifold as in \cite{barnettquadratic} with a manually determined regularization coefficient $\lambda = 0.5$ ($\alpha$ in (27) of \cite{barnettquadratic}) and the rational manifold using the fully nonlinear curve-fitting approach outlined in the previous section. In \autoref{fig:burgersreconstr} we display the reconstruction in addition to the original data using an $x-t$ plot. Furthermore, in \autoref{fig:burgersreconstrerror} we plot the local error in space--time defined, with some abuse of notation, as:
\begin{equation*}
    \varepsilon_{xt} = \bm{\varphi}(A) - X.
\end{equation*}

It can be seen that for $r=15$, the reconstruction accuracy of both the quadratic as POD linear manifold is poor, whereas the reconstruction of the rational quadratic manifold is visually nearly identical to the data. Indeed, the largest local error of the rational quadratic manifold is at most approximately $0.003$ which is two orders of magnitude lower than the largest errors of the quadratic and linear manifolds (approximately $0.4$ and $0.6$ respectively). The sources of error for the linear and quadratic manifolds are predominantly oscillations around the moving shock discontinuity as can be seen in \autoref{fig:burgersreconstrerror}. This shows the poor performance of these methods for such problems. The rational quadratic manifold also oscillates around the shock, but with a much smaller amplitude, indicating that it is better-suited for shock-dominated problems. The accuracy of the rational quadratic manifold is much higher than the quadratic and POD linear manifolds. 

\begin{figure}
    \centering
    \includegraphics[width = \textwidth]{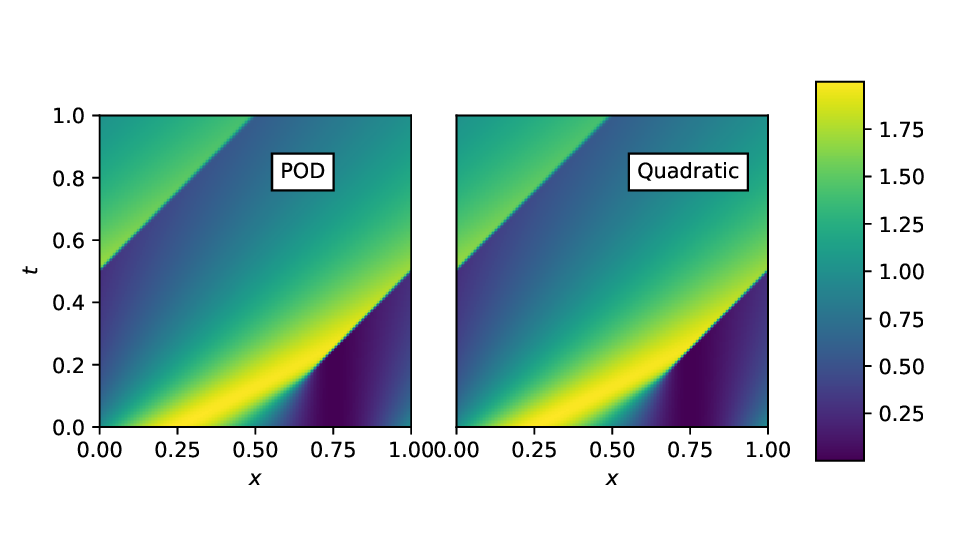}
    \caption{Reconstructions for $\varepsilon_{xt}^{\text{max}} \approx 3\cdot 10^{-3}$ in space-time plots.}
    \label{fig:burgersreconstr2}
\end{figure}

We note that the increased accuracy comes at the cost of a slow fitting procedure. The precise fitting times and maximum absolute space-time errors, $\varepsilon_{xt}^{\text{max}} := \max_{i,j}|\varepsilon_{xt}|_{i,j}$, have been displayed in \autoref{tab:burgersperf}.
\begin{table}[]
\centering
\begin{tabular}{llll}
\toprule
            & $t_{\text{fit}}$ & $t_{\text{online}}$ & $\varepsilon_{xt}^{\text{max}}$ \\ \midrule
entropy stable rational quadratic ($r=15) $ & 1851 s   &  7.01 s  & $3.36 \cdot 10^{-3}$                       \\
rational quadratic ($r=15) $ & 1851 s  &  6.74 s  & $3.36 \cdot 10^{-3}$                       \\
POD linear ($r=15$)    & 0.107 s  & 0.543 s & 0.675                                                                      \\
POD linear ($r=160$)    & 0.107 s &  0.622 s & $2.61\cdot 10^{-3}$                                                                     \\
quadratic ($r=15$)      & 0.964 s & 6.14 s & 0.456                                                                      \\
quadratic ($r=150$)      & 122 s  & 121 s & $4.70 \cdot 10^{-3}$                                                                     \\ \bottomrule
\end{tabular}
\caption{Accuracy and computational cost for fitting Burgers solution with different manifold parameterizations}
\label{tab:burgersperf}
\end{table}
When we construct the POD linear manifold and quadratic manifold to an accuracy of $\varepsilon_{xt}^{\text{max}} \approx 3 \cdot 10^{-3}$ we see that we require approximately $r = 160$ and $r = 150$, respectively. The reconstructions are given in \autoref{fig:burgersreconstr2} and the largest space--time errors $\varepsilon_{xt}$ are $\mathcal{O}(4\cdot 10^{-3})$. The changes in fitting time have also been denoted in \autoref{tab:burgersperf}. A large increase can be observed for quadratic manifolds, while the linear POD stays constant as the implementation calculates all $n_s$ singular vectors at once.

\begin{figure}[h!]
    \centering
    \begin{minipage}{0.5\textwidth}
        \centering
        \includegraphics[width = \textwidth]{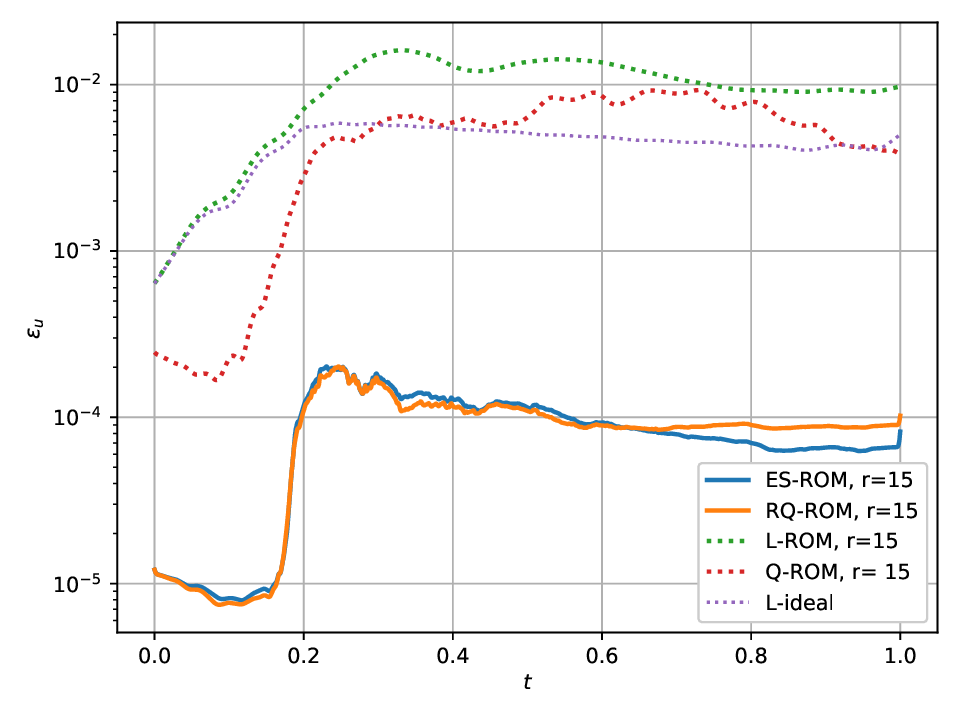}
        \caption{Total error evolution for the ROMs $\varepsilon_u(t)$ and the ideal linear projection error $\varepsilon_{\text{proj}}(t)$ for reduced manifolds with constant $r = 15$.}
        \label{fig:burgerserrorevor}
    \end{minipage}\hfill
    \begin{minipage}{0.5\textwidth}
        \centering
        \includegraphics[width = \textwidth]{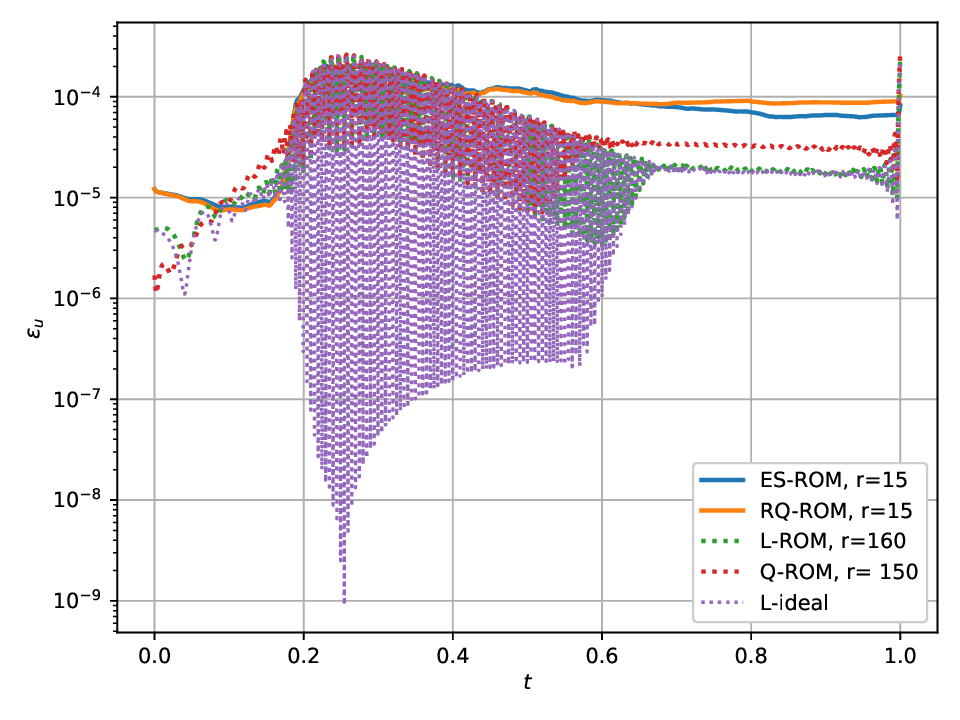}
        \caption{Total error evolution for the ROMs $\varepsilon_u(t)$ and the ideal linear projection error $\varepsilon_{\text{proj}}(t)$ for reduced manifolds with approximately constant $\varepsilon_{xt}^{\text{max}}$.}
        \label{fig:burgerserrorevoeps}
    \end{minipage}
\end{figure}

\begin{figure}[h!]
    \centering
    \begin{minipage}{0.5\textwidth}
        \centering
        \includegraphics[width = \textwidth]{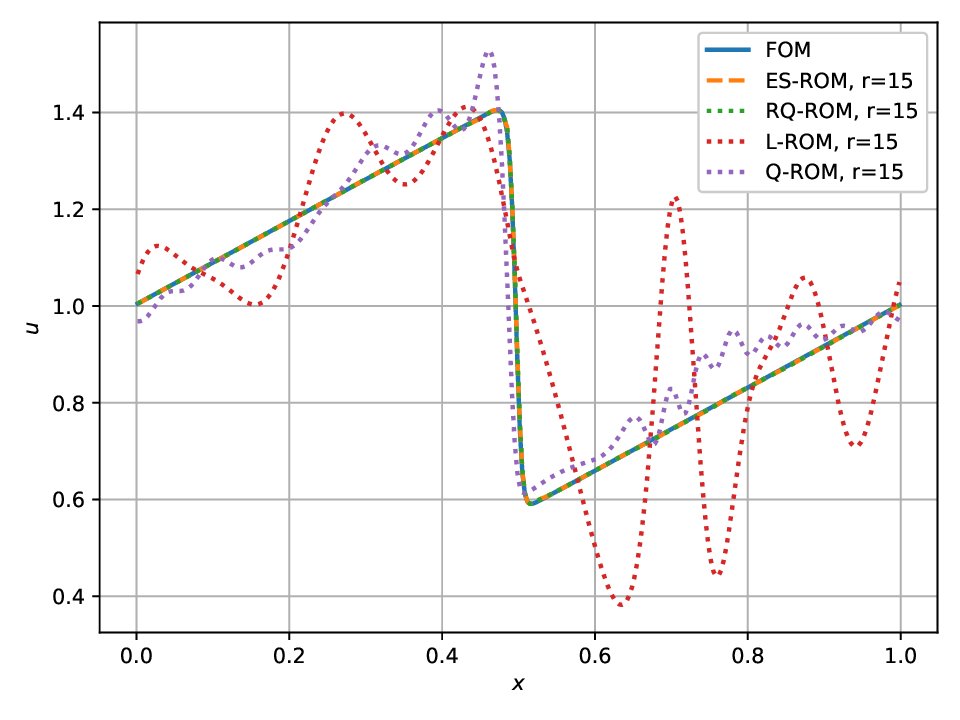}
        \caption{Spatial profile $\bm{u}_r(t)$ for the ROMs at $t = T$ for reduced manifolds with constant $r = 15$.}
        \label{fig:burgersprofiler}
    \end{minipage}\hfill
    \begin{minipage}{0.5\textwidth}
        \centering
        \includegraphics[width = \textwidth]{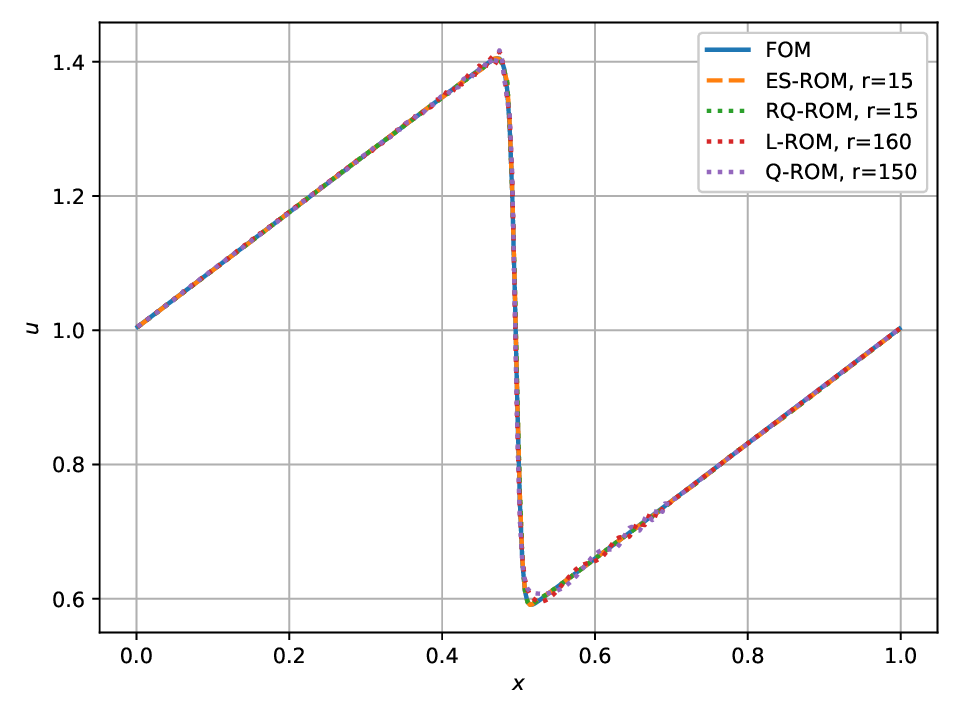}
        \caption{Spatial profile $\bm{u}_r(t)$ for the ROMs at $t = T$ for reduced manifolds with approximately constant $\varepsilon_{xt}^{\text{max}}$.}
        \label{fig:burgersprofileeps}
    \end{minipage}
\end{figure}

We continue to consider the ROM performance and accuracy in more detail. We compare the rational quadratic manifold ROM in entropy stable \eqref{eq:esROM} (ES-ROM) and generic \eqref{eq:galrom} (RQ-ROM) form to a linear manifold POD-Galerkin ROM (L-ROM) and a quadratic manifold Galerkin ROM (Q-ROM). We will make one comparison of the ROMs using the previously obtained manifolds with $r = 15$ and another comparison where $r$ is chosen such that each manifold has approximately the same space--time reconstruction error $\varepsilon_{xt}^{\text{max}} \approx 3 \cdot 10^{-3}$. The initial conditions for the simulations with $r = 15$ will be taken as the first column of the matrix $A = \Phi^TX$ and the entropy stable form of the rational quadratic ROM will have $\alpha = 0$ at $t = 0$. For the manifolds that have approximately the same accuracy we will take the first columns of the matrices $A_r = \Phi_r^T X \in \mathbb{R}^{r \times n_s}$ with $r = 160$ and $r = 150$ for the linear and quadratic manifolds, respectively. Here, $\Phi_r \in \mathbb{R}^{N_h \times r}$ are the $r$ first singular vectors of $X$. We plot the temporal evolution of the total error norm:
\begin{equation*}
    \varepsilon_u(t) := ||\bm{u}_h(t) - \bm{u}_r(t)||_{\Omega_h},
\end{equation*}
and the ideal linear projection error (L-ideal):
\begin{equation*}
    \varepsilon_\text{proj}(t) := ||(I - \Pi_{r_{\text{lin}}}) \bm{u}_h(t)||_{\Omega_h},
\end{equation*}
where $\Pi_{r_{\text{lin}}} : \mathbb{R}^{N_h} \rightarrow \mathcal{V}$ projects on the $r_{\text{lin}}$-dimensional reduced space of the respective linear ROMs in the different experiments. The ideal projection error forms a lower bound for the linear POD-Galerkin ROM error. To measure computational performance we will track the runtimes of the online phases which we denote $t_{\text{online}}$. The results for the simulations with constant $r$ are shown in \autoref{fig:burgerserrorevor} and the results for approximately constant space--time error are given in \autoref{fig:burgerserrorevoeps}. We also show the spatial profile of the solution at $t = T$ as predicted by the different ROMs and different manifolds in \autoref{fig:burgersprofiler} for constant $r$ and \autoref{fig:burgersprofileeps} for constant space--time error, respectively. The rational manifold ROMs clearly outperform the others in case of constant $r$ and the differences between the results of the rational manifolds themselves are nearly zero. Steep increases in errors occur for all ROMs upon formation of the shock which indicates this is indeed a difficult instant of the flow for the reduced manifolds to fit to the data. For the simulation with manifolds with constant error the performance of the ROMs is more equal, with the linear and quadratic ROMs only suffering of some oscillations in the spatial profile. The oscillations occur at moments when the ideal projections error $\varepsilon_{\text{proj}}$ is also oscillatory in time. The rational quadratic manifolds do not suffer from oscillations. Because of the large reduced spaces required to obtain approximately the same reconstruction errors as the rational manifolds, the linear and quadratic manifold ROMs have more expensive online phases than when tested at constant $r$. This is especially notable for quadratic manifolds where computing the Jacobian and its Moore-Penrose pseudoinverse contribute heavily to the increase in cost. At constant $r$ the quadratic manifold ROM and the rational manifold ROMs are nearly equally fast showing that the Jacobian of the quadratic manifold parameterization and of the rational manifolds with and without enrichment are nearly equally expensive to compute. 
From this experiment we conclude that at the cost of a slower fitting process rational quadratic manifolds can significantly outperform quadratic and linear POD manifolds in terms of reconstruction accuracy for the same number reduced space dimensions. 

\subsection{Shallow water equations}\label{sec:shallow}
This experiment will mainly focus on the entropy conservation and stability aspect of our proposed ROM on rational manifolds. We show that our novel entropy stable ROM satisfies the reduced total entropy estimate \eqref{eq:romesstability}. To this end we will carry out experiments with the shallow water equations. We use the shallow water equations due to their nontrivial entropy function, which will be defined later, as compared to the Burgers equation. We will perform one experiment where discontinuities appear in the solution and one where the solution remains smooth during the time interval of interest. In the smooth case we can run the FOM and ROMs without entropy dissipation operators. As a result we can analyse the entropy conservation properties of the ROM. In the discontinuous case we will analyse the behaviour of the reduced entropy as compared to the FOM entropy.

In the following we briefly introduce the shallow water equations and the entropy stable numerical scheme used to obtain the FOM. The shallow water equations are given by:
\begin{equation}
    \frac{\partial}{\partial t}\begin{bmatrix}
        h \\ hu
    \end{bmatrix} + \frac{\partial}{\partial x}\begin{bmatrix}
        hu \\ hu^2 + \frac{1}{2}gh^2
    \end{bmatrix} = 0,
    \label{eq:shallow}
\end{equation}
with conserved variables $h : \Omega \times [0,T] \rightarrow \mathbb{R}$ and $hu : \Omega \times [0,T] \rightarrow \mathbb{R}$ representing the local water column height and the momentum per unit mass, respectively. We collect the conserved variables in a vector $\bm{u} := \left[h, hu\right]^T$. The constant $g \in \mathbb{R}^+$ is the positive real gravitational acceleration, we assume the equations are normalized such that $g = 3$, which gave challenging test cases for our spatial domain size and initial conditions. As continuous and discrete total entropy we take the common choice \cite{fjordholmenergy}:
\begin{equation*}
    \mathcal{S}[\bm{u}] = \int_{\Omega}\frac{1}{2}\left(\frac{u_2^2}{u_1} + gu_1^2\right)dx,
\end{equation*}
leading to:
\begin{equation*}
    S_h[\bm{u}_h] = \sum_i \Delta x_i \frac{h_i u_i^2 + gh_i^2}{2},
\end{equation*}
from which the reduced total entropy follows. This choice leads to the following entropy variables:
\begin{equation*}
    \bm{\eta}(\bm{u}) = \begin{bmatrix}
        g u_1 - \frac{1}{2}\left(\frac{u_2}{u_1}\right)^2 \\ \frac{u_2}{u_1}
    \end{bmatrix}, 
\end{equation*}
where $u_1 = h$ and $u_2 = hu$, and inverse function:
\begin{equation*}
    \bm{u}(\bm{\eta}) =\frac{2\eta_1 + \eta_2^2}{2g} \begin{bmatrix}
        1 \\ \eta_2
    \end{bmatrix}.
\end{equation*}
An entropy conservative flux is given by \cite{fjordholmenergy, fjordholmarbitrarily}:
\begin{equation*}
    \bm{f}_{i+1/2} = \begin{bmatrix}
        \overline{h}_{i+1/2}\overline{u}_{i+1/2} \\ \overline{h}_{i+1/2}\cdot \overline{u}_{i+1/2}^2 + \frac{1}{2}g\overline{h^2}_{i+1/2}
    \end{bmatrix},
\end{equation*}
where $\overline{a}_{i+1/2} = \frac{1}{2}(a_{i+1} + a_i)$ indicates taking the average of neighbouring volume based quantities. As entropy dissipation operator we take the diffusion operators of Roe type (see \cite{fjordholmarbitrarily}) with the eigenvalues and eigenvectors of the flux Jacobian evaluated at the arithmetic average of neighbouring values. We obtain a second accurate entropy dissipation operator using the entropy stable total variation diminishing (TVD) reconstruction based on the minmod limiter (see \cite{fjordholmarbitrarily}).

We will discretize \eqref{eq:shallow} for both experiments on a domain $\Omega = \mathbb{T}([-L,L])$ with $L = 1$ using a numerical grid consisting of 300 cells. Discretization in time will be done using the RK4 method with a time step of size $\Delta t = 0.0005$. We will integrate the discontinuous experiment in time until $T = 1$ and the smooth experiment until $T = 0.5$. Solution snapshots are captured every 5 timesteps resulting in $n_s = 401$ and $n_s = 201$ for the discontinuous and smooth experiment, respectively. For the discontinuous case we will be interested in a dam break problem, this means we will take as initial condition:
\begin{equation*}
    h_0(x) = \begin{cases}
        1.5 & |x|<0.2, \\
        1 & |x|\geq 0.2,
    \end{cases} \quad (hu)_0(x) = 0.
\end{equation*}
The smooth case will consist of a quiescent water level with a small perturbation, so that the initial condition is given by:
\begin{equation*}
    h_0(x) = 1 + 0.1 \cdot \exp\left(-100 \cdot x^2\right), \quad (hu)_0(x) = 0.
\end{equation*}
The reduced space dimension is taken at $r = 15$ for both experiments. 

\begin{figure}
    \centering
    \includegraphics[width = \textwidth]{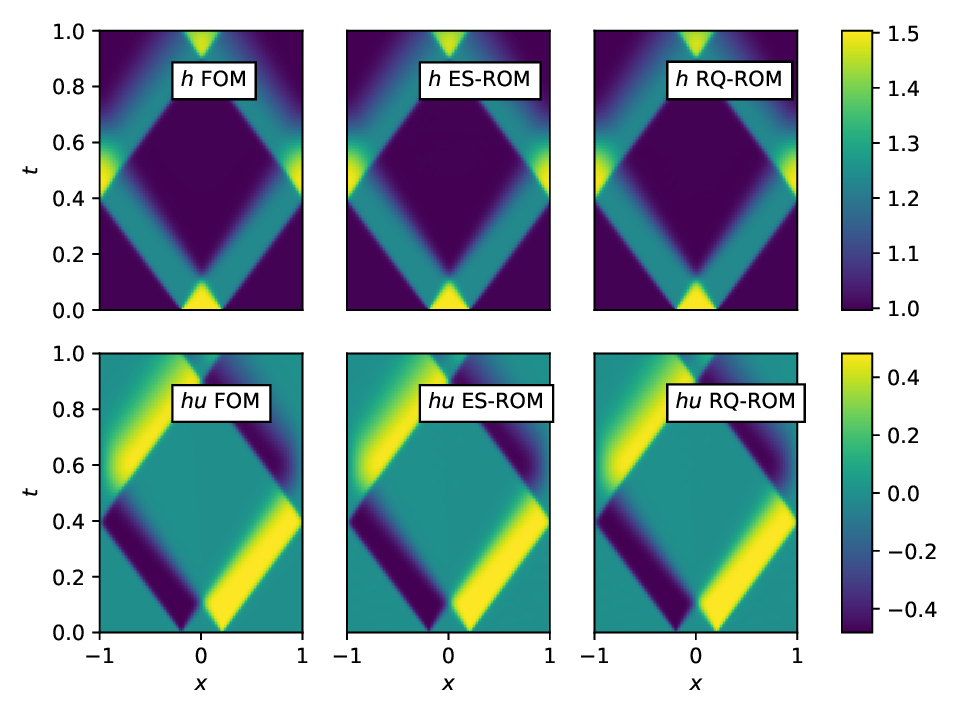}
    \caption{$x-t$ plot of solution approximation by the ROMs and FOM for the dam break problem.}
    \label{fig:shallowxt1}
\end{figure}

\begin{figure}
    \centering
    \includegraphics[width = \textwidth]{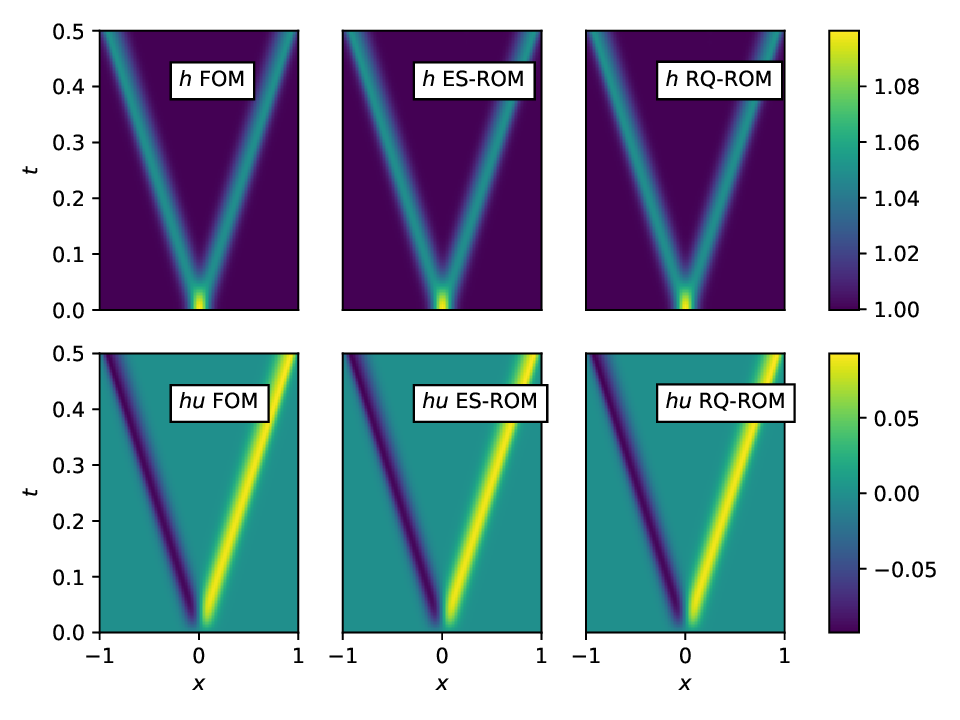}
    \caption{$x-t$ plot of solution approximation by the ROMs and FOM for the water height perturbation problem.}
    \label{fig:shallowxt2}
\end{figure}

\begin{figure}[h!]
    \centering
    \begin{minipage}{0.5\textwidth}
        \centering
        \includegraphics[width = \textwidth]{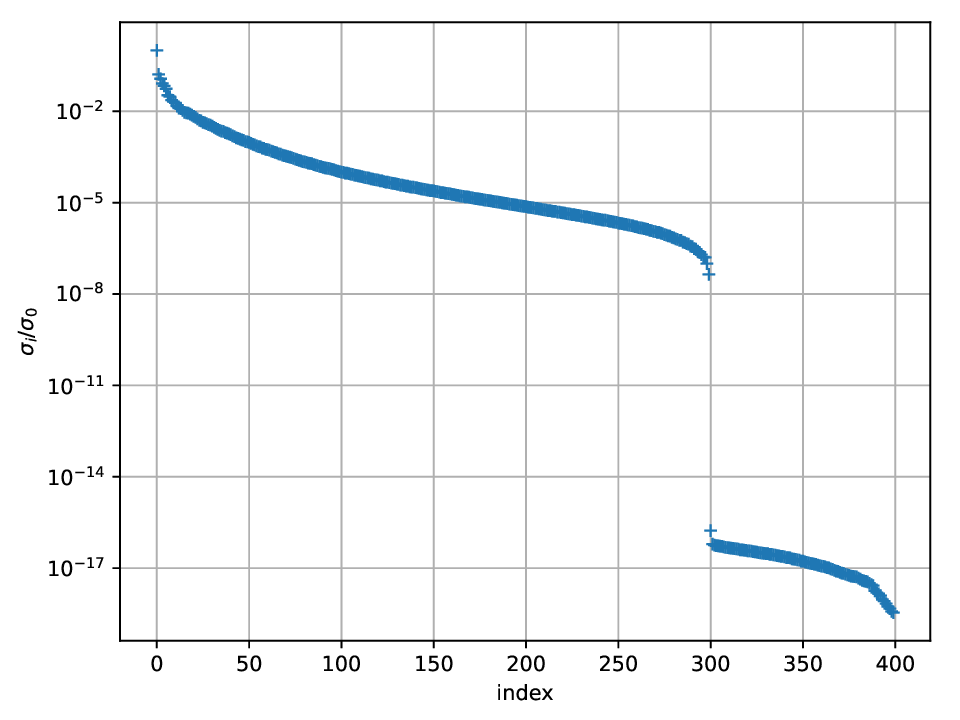}
        \caption{Normalized singular values of the shallow water equations data for the dam break problem.}
        \label{fig:shallowsingvals1}
        \end{minipage}\hfill
    \begin{minipage}{0.5\textwidth}
        \centering
        \includegraphics[width = \textwidth]{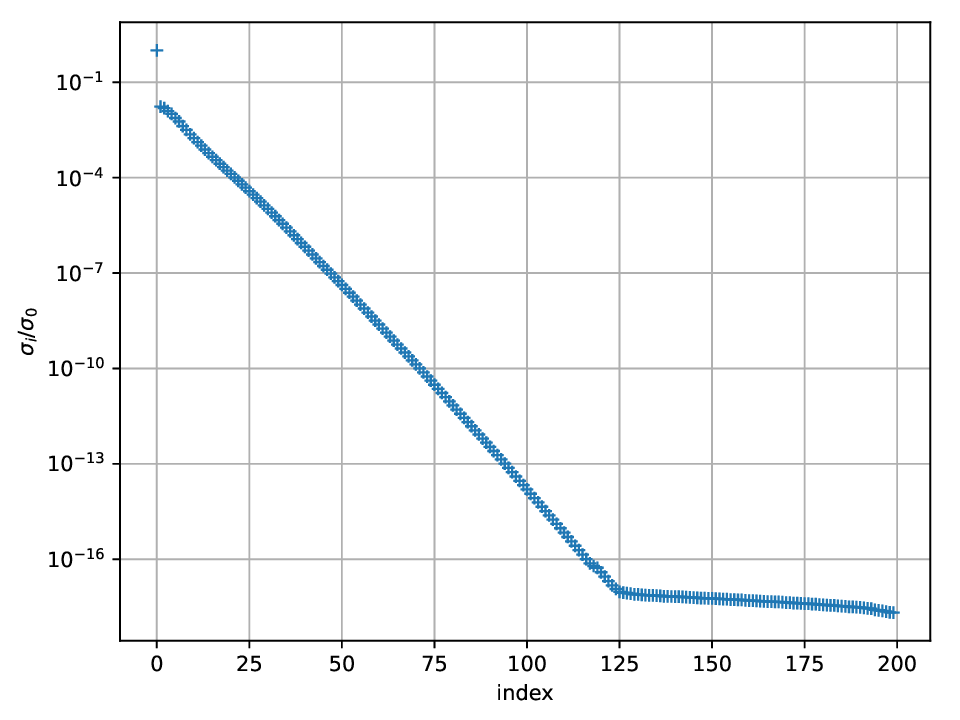}
        \caption{Normalized singular values of the shallow water equations data for the water height perturbation problem.}
        \label{fig:shallowsingvals2}
    \end{minipage}
\end{figure}

Although it is important that entropy is appropriately conserved or dissipated, accuracy with respect to the FOM is also required for an effective ROM. Hence, before we analyse the conservation properties of our proposed ROM, we compare the FOM solution approximation quality of our entropy stable ROM and the generic ROM. We will provide space--time plots of both the discontinuous and smooth experiments. The discontinuous case is given in \autoref{fig:shallowxt1} and the smooth case is given in \autoref{fig:shallowxt2}. Visually, both ROMs closely resemble the FOM in both cases. Furthermore, it can be seen from the sharp color gradients that the dam break problem develops shocks. To emphasize these are difficult cases for linear model reduction approaches we also plot the normalized singular value decay in \autoref{fig:shallowsingvals1} and \autoref{fig:shallowsingvals2} for the dam break and water height perturbation problems, respectively. For the dam beak problem the decay is very slow and the water height perturbation decays moderately slow, indicating slow and moderately slow Kolmogorov $n$-width decay \eqref{eq:knW}.

\begin{figure}[h!]
    \centering
    \includegraphics[width = \textwidth]{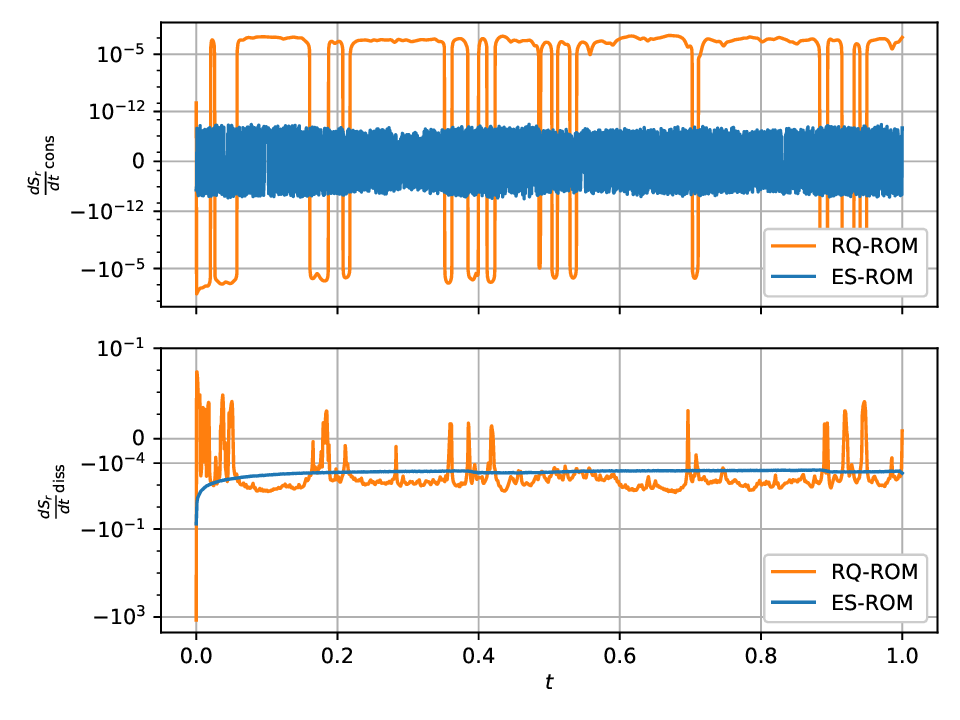}
    \caption{The entropy production $\left(\frac{dS_r}{dt}\right)_{\text{cons}}$ of the conservative part of the ROMs and the entropy production $\left(\frac{dS_r}{dt}\right)_{\text{diss}}$ of the dissipative part of the ROMs for the dam break problem.}
    \label{fig:shallowentropyprod}
\end{figure}

\begin{figure}[h!]
    \centering
    \includegraphics[width = \textwidth]{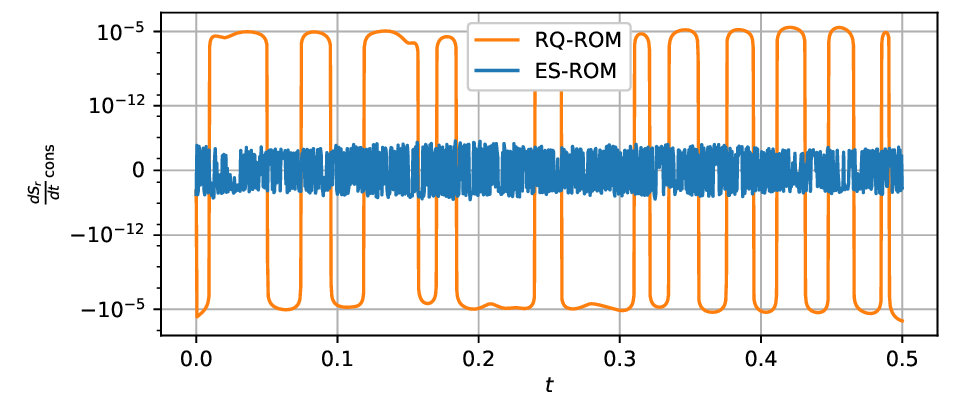}
    \caption{The entropy production $\left(\frac{dS_r}{dt}\right)_{\text{cons}}$ of the conservative part of the ROMs for the water height perturbation problem.}
    \label{fig:shallowentropyprod2}
\end{figure}

We will analyse the entropy conservation and stability properties of the entropy stable (ES-ROM) and generic (RQ-ROM) rational manifold ROMs. To this end, we define the entropy error:
\begin{equation*}
    \varepsilon_{\mathcal{S}}(t) := \left|S_h[\bm{u}_h(t)] - S_r[\bm{a}(t)] \right|,
\end{equation*}
giving the absolute instantaneous deviation of the entropy of the ROM from the entropy of the FOM. Similarly we will define the entropy conservation error:
\begin{equation*}
    \varepsilon_{\mathcal{S}_0}(t) := \left|S_r[\bm{a}(0)] - S_r[\bm{a}(t)] \right|,
\end{equation*}
which measures the departure from the initial entropy and thus the error in exact conservation of the entropy in time. Since our models are semi-discretely entropy stable we have to monitor the instantaneous time rate of change of reduced total entropy \eqref{eq:romesstability} to verify that our proposed theoretical framework works. Hence, we will analyse the contribution to the total entropy production \eqref{eq:entropyprod} of two separate parts of the ROM \eqref{eq:esROM}, namely the entropy conserving part:
\begin{equation*}
    \left(\frac{dS_r}{dt}\right)_{\text{cons}} := -\tilde{\bm{\eta}}_r^T\Delta_v \bm{f}_h^*(\bm{u}(\tilde{\bm{\eta}}_r)),
\end{equation*}
which should equal zero to machine precision, and the entropy dissipative part:
\begin{equation*}
    \left(\frac{dS_r}{dt}\right)_{\text{diss}} := \tilde{\bm{\eta}}_r^T\Delta_vD_h(\tilde{\bm{u}}_r)\Delta_i\tilde{\bm{\eta}}_r,
\end{equation*}
which should always be negative or zero. Similar quantities can be defined in an obvious manner for the generic ROM without entropy projection. In the results given by \autoref{fig:shallowentropyprod} and \autoref{fig:shallowentropyprod2} we have used symmetric log plots which are linear around zero so that negative values can also be plotted. This allows us to see when a ROM is unphysically producing entropy i.e.\ $\left(\frac{dS_r}{dt}\right)_{\text{cons}}, \left(\frac{dS_r}{dt}\right)_{\text{diss}} > 0$.

\begin{figure}[h!]
    \centering
    \includegraphics[width = \textwidth]{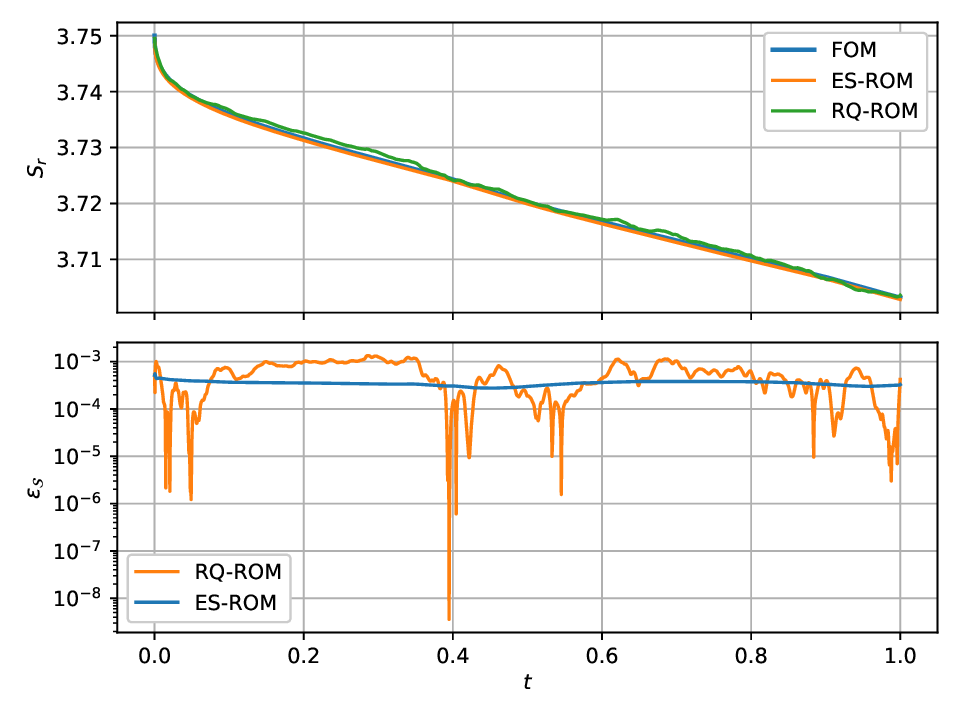}
    \caption{The evolution of the reduced total entropy $\mathcal{S}_r$ of the entropy stable and generic ROM compared to the entropy of the FOM and the evolution of the entropy error $\varepsilon_{\mathcal{S}}$ of the entropy stable and generic ROM for the dam break problem.}
    \label{fig:shallowentropyevo}
\end{figure}

The results of the discontinuous dam break experiment are displayed in \autoref{fig:shallowentropyevo} and \autoref{fig:shallowentropyprod}. The results confirm that the proposed entropy stable framework works as expected. This is the case since the entropy production of the conservative part is zero to machine precision for the entropy stable ROM while the entropy dissipative part does not change sign and is indeed negative. The entropy production of the conservative part of the generic ROM is orders of magnitude larger than that of the entropy stable ROM. An important point is that large portions in time of the entropy production by the entropy conservative part are positive (instead of zero). This indicates physically incorrect behaviour as entropy is being produced. The contribution of the entropy dissipation operator from the generic ROM is erratic. Moreover, it is also occasionally positive showing that this part of the ROM is also sometimes producing physically incorrect results. The temporal evolution of the entropy error, $\varepsilon_{\mathcal{S}}$, is also given in \autoref{fig:shallowentropyevo}. The temporal evolution is approximately constant in time for the entropy stable ROM and $\varepsilon_\mathcal{S}$ is small. This indicates that the evolution of the entropy behaves roughly the same as the FOM and is off mainly due to an error in representation of the initial condition and of subsequent FOM solutions. The general behaviour of the entropy error of the generic ROM is erratic and shows that the entropy of the generic ROM oscillates around the values predicted by the FOM. This can also be seen in the temporal evolution of the reduced entropy as in the top panel of \autoref{fig:shallowentropyevo}.

\begin{figure}[h!]
    \centering
    \includegraphics[width = \textwidth]{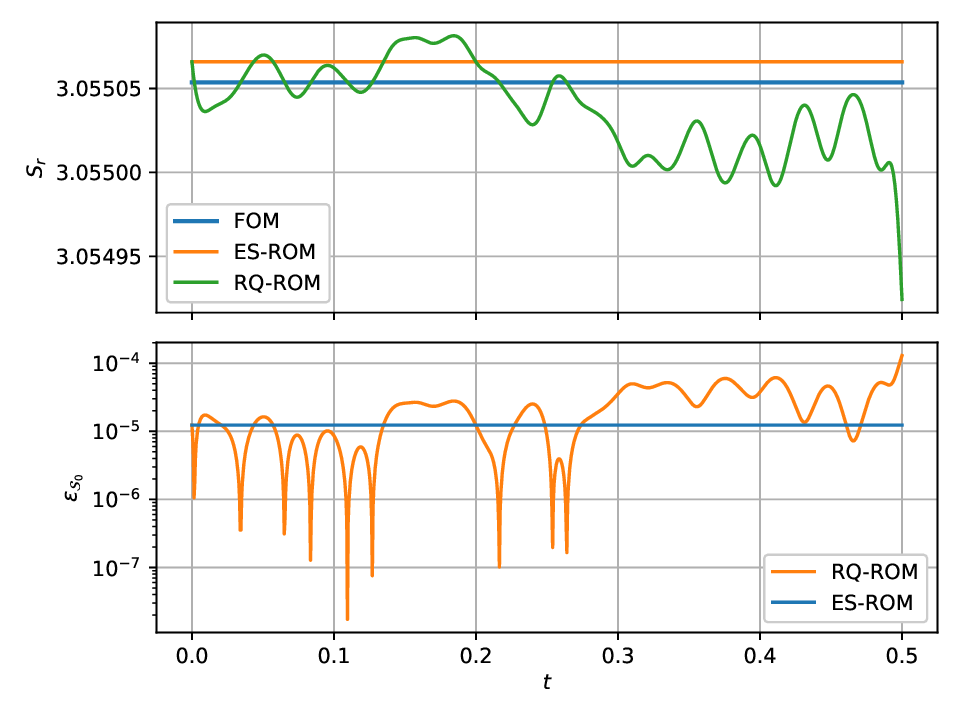}
    \caption{The evolution of the reduced total entropy $\mathcal{S}_r$ of the entropy conserving and generic ROM compared to the entropy of the FOM and the evolution of the entropy error $\varepsilon_{\mathcal{S}_0}$ of the entropy conserving and generic ROM for the water height perturbation problem.}
    \label{fig:shallowentropyevo2}
\end{figure}

The results of the smooth experiment are shown in \autoref{fig:shallowentropyevo2} and \autoref{fig:shallowentropyprod2}. As there is no entropy dissipation present in the FOM or ROM the entropy should remain approximately constant (exact conservation is difficult since RK4 is not an entropy conservative time-integrator \cite{lozanoentropy, lozanoentropy2, tadmorentropy}). For our entropy stable ROM this is indeed the case as can be seen from the bottom panel of \autoref{fig:shallowentropyevo2}, the entropy conservation error, $\varepsilon_{\mathcal{S}_0}$, does not exceed $\mathcal{O}(10^{-9})$. The entropy of our entropy stable ROM stays almost exactly constant. The error in entropy with respect to the FOM is almost entirely dictated by the initial error $\varepsilon_{\mathcal{S}}(0) = \left|S_h[\bm{u}_h(t)] - S_r[\bm{a}(0)]\right|$. The generic ROM does not conserve entropy and its entropy conservation error behaves erratically. This manifests itself in clear deviations from the FOM entropy which can be observed in the top panel of \autoref{fig:shallowentropyevo2}. For completeness we also plot the entropy production of the conservative part of the spatial discretization of the ROM $\left(\frac{dS_r}{dt}\right)_{\text{cons}}$ in \autoref{fig:shallowentropyprod2}. Again, it can be seen that entropy is conserved up to machine precision by the spatial discretization of our entropy stable ROM where this is not the case for the generic ROM. Additionally, the generic ROM produces entropy during several intervals of the simulation and is therefore not physically correct. From both experiments, we conclude that our novel entropy stable ROM ensures physically correct behaviour, whereas this cannot be assumed for the generic manifold Galerkin ROM.

\subsection{Compressible Euler equations}\label{sec:euler}
The focus of this experiment is on the effect of the entropy projection and tangent space enrichment on the accuracy of the ROM. We will be interested in particular in the benefit of tangent space enrichment in the reconstruction accuracy of the entropy projection. In addition, we will analyse the error that can be incurred with respect to the FOM by the introduction of an entropy projection step in the ROM as we propose. A good case to study for this experiment are the compressible Euler equations. Due to their nontrivial entropy functional and corresponding entropy variables it is not expected that without extra measures, like TSE, the entropy projection will be accurate.

A short introduction to the compressible Euler equations now follows. The compressible Euler equations are given by:
\begin{equation}
    \frac{\partial}{\partial t}\begin{bmatrix}
        \rho \\ \rho u \\ E
    \end{bmatrix} + \frac{\partial}{\partial x}\begin{bmatrix}
        \rho u \\ \rho u^2 + p \\ (E + p)u
    \end{bmatrix} = 0,
    \label{eq:euler}
\end{equation}
where $\rho : \Omega \times [0,T] \rightarrow \mathbb{R}$ is the density, $\rho u : \Omega \times [0,T] \rightarrow \mathbb{R}$ is the momentum and $E : \Omega \times [0,T] \rightarrow \mathbb{R}$ is the total energy. We gather these conserved variables in a vector $\bm{u} = \left[\rho, \rho u, E\right]^T$. Furthermore, we assume the equations are suitably normalized so that they are dimensionless. The pressure $p : \mathbb{R}^n \rightarrow \mathbb{R}$ is related to the conserved variables through an equation of state, representing the thermodynamics at hand:
\begin{equation}
    p(\bm{u}) = \left(u_3 - \frac{1}{2}\frac{u_2^2}{u_1}\right)(\gamma - 1),
    \label{eq:eos}
\end{equation}
where $\gamma \in \mathbb{R}^+$ is the specific heat ratio, which we take at the standard choice $\gamma = 1.4$. The thermodynamic quantities, i.e.\ pressure, density and total energy are necessarily nonnegative. Throughout the experiments we will assume our FOM and ROMs respect this condition, assuring this condition mathematically may be the subject of future work. The entropy functional of interest will be taken as:
\begin{equation*}
    \mathcal{S}[\bm{u}] = \int_{\Omega} \frac{-u_1 \sigma}{\gamma - 1} dx,
\end{equation*}
where $\sigma : \mathbb{R}^n \rightarrow \mathbb{R}$ is the specific entropy defined as a function of the conserved variables like:
\begin{equation*}
    \sigma(\bm{u}) = \ln \left(\frac{p}{u_1^\gamma}\right) ,
\end{equation*}
where $p$ is evaluated using \eqref{eq:eos}. Different entropies are also possible, see for instance \cite{hartensymmetric}. In turn, this gives rise to the discrete total entropy functional:
\begin{equation*}
    S_h[\bm{u}_h] = \sum_i \Delta x_i \frac{-\rho_i \sigma_i}{\gamma - 1}.
\end{equation*}
The associated entropy variables are given by:
\begin{equation*}
    \bm{\eta}(\bm{u}) = \begin{bmatrix}
        \frac{\gamma - \sigma}{\gamma - 1} - \frac{u_2^2}{2u_1 p} \\ u_2/p \\ -u_1/p
    \end{bmatrix},
\end{equation*}
and consequently the inverse of the entropy variables is:
\begin{equation*}
    \bm{u}(\bm{\eta}) =  \exp\left(\frac{\gamma}{1-\gamma} - \left[\eta_1 - \frac{1}{2}\frac{\eta_2^2}{\eta_3} \right] \right)\begin{bmatrix}
       (-\eta_3)^{\frac{1}{1-\gamma}} \\ -\eta_2(-\eta_3)^{\frac{\gamma}{1-\gamma}} \\ \left(\frac{1}{2}\eta_2^2 - \frac{\eta_3}{\gamma - 1} \right)\cdot (-\eta_3)^{\frac{2\gamma - 1}{1-\gamma}}
    \end{bmatrix}.
\end{equation*}
Considering $-\eta_3 = \rho / p$ is generally exponentiated to some noninteger power we see the importance of positivity of the thermodynamic variables. As in \cite{fjordholmarbitrarily}, we will use the entropy conserving numerical flux by Ismail and Roe \cite{ismailaffordable} for which we define the following variables:
\begin{equation*}
    \bm{z} = \begin{bmatrix}
        z^1 \\ z^2 \\ z^3
    \end{bmatrix} = \sqrt{\frac{\rho}{p}}\begin{bmatrix}
        1 \\ u \\ p
    \end{bmatrix},
\end{equation*}
finally an entropy conservative flux is given by $\bm{f}_{i+1/2} = \begin{bmatrix} \bm{F}_{i+1/2}^1 & \bm{F}_{i+1/2}^2 & \bm{F}_{i+1/2}^3 \end{bmatrix}^T$:
\begin{align*}
    \bm{F}_{i+1/2}^1 &= \overline{z^2}_{i+1/2} \cdot (z^3)_{i+1/2}^{\ln}, \\ 
    \bm{F}_{i+1/2}^2 &= \frac{\overline{z^3}_{i+1/2}}{\overline{z^1}_{i+1/2}} + \frac{\overline{z^2}_{i+1/2}}{\overline{z^1}_{i+1/2}}\bm{F}_{i+1/2}^1, \\
    \bm{F}_{i+1/2}^3 &= \frac{1}{2}\frac{\overline{z^2}_{i+1/2}}{\overline{z^1}_{i+1/2}}\left(\frac{\gamma + 1}{\gamma - 1}\frac{(z^3)_{i+1/2}^{\ln}}{(z^1)_{i+1/2}^{\ln}} + \bm{F}_{i+1/2}^2\right),
\end{align*}
where $a^{\ln}$ denotes the logarithmic mean, which is defined as:
\begin{equation*}
    a_{i+1/2}^{\ln} := \frac{a_{i+1} - a_{i}}{\ln{a_{i+1}} - \ln{a_i}}.
\end{equation*}
Computation of the logarithmic mean is generally not numerically stable when $a_{i+1}\approx a_i$, but a popular algorithm which we will use to deal with this is also given in \cite{ismailaffordable}. There is an abundance of alternative entropy conservative numerical fluxes that can be used \cite{ranochacomparison, chandrashekarkinetic, hickenfamily, tadmornumerical, kuyahigh, tamakicomprehensive} some of which also conserve kinetic energy in the sense of Jameson \cite{jamesonformulation}. As an entropy dissipation operator we take the Roe-type diffusion operator \cite{fjordholmarbitrarily} where the eigenvalues and vectors of the flux Jacobian are evaluated at the arithmetic average of the neighbouring conserved variables. As with the shallow water equations, we obtain a second accurate entropy dissipation operator using the entropy stable total variation diminishing (TVD) reconstruction based on the minmod limiter \cite{fjordholmarbitrarily}.

For the experiment we will consider a periodic modification of the famous Sod's shock tube \cite{sodsurvey}, which avoids the need to implement entropy stable boundary conditions. We will discretize \eqref{eq:euler} on a domain $\Omega = \mathbb{T}([0,L])$ with $L = 1$ on a numerical grid of $250$ cells. The number of cells is relatively small to facilitate a relatively short manifold learning process. Integration of the ROM in time will be carried out using the RK4 time integrator with a time step size $\Delta t = 0.0001$. We will integrate in time until $T = 0.5$ (beyond the typical time used for this experiment), resulting in interesting shock-rarefaction interactions. Again, we will capture snapshots after every 5 timesteps so that we have $n_s = 1001$. Our periodic modification of Sod's shock tube experiment has an initial condition given by:
\begin{equation*}
    \rho_0(x) = \begin{cases}
        1 & 0.25 < x < 0.75 \\
        0.125 & \text{elsewhere}
    \end{cases} \quad u_0(x) = 0, \quad p_0(x) = \begin{cases}
        1 & 0.25 < x < 0.75 \\
        0.1 & \text{elsewhere}
    \end{cases},
\end{equation*}
where the conserved variables $(\rho, \rho u, E)$ are calculated from these primitive variables using the equation of state \eqref{eq:eos} and the definition of momentum. We take a reduced space dimension $r = 15$.

\begin{figure}
    \centering
    \includegraphics[width = \textwidth]{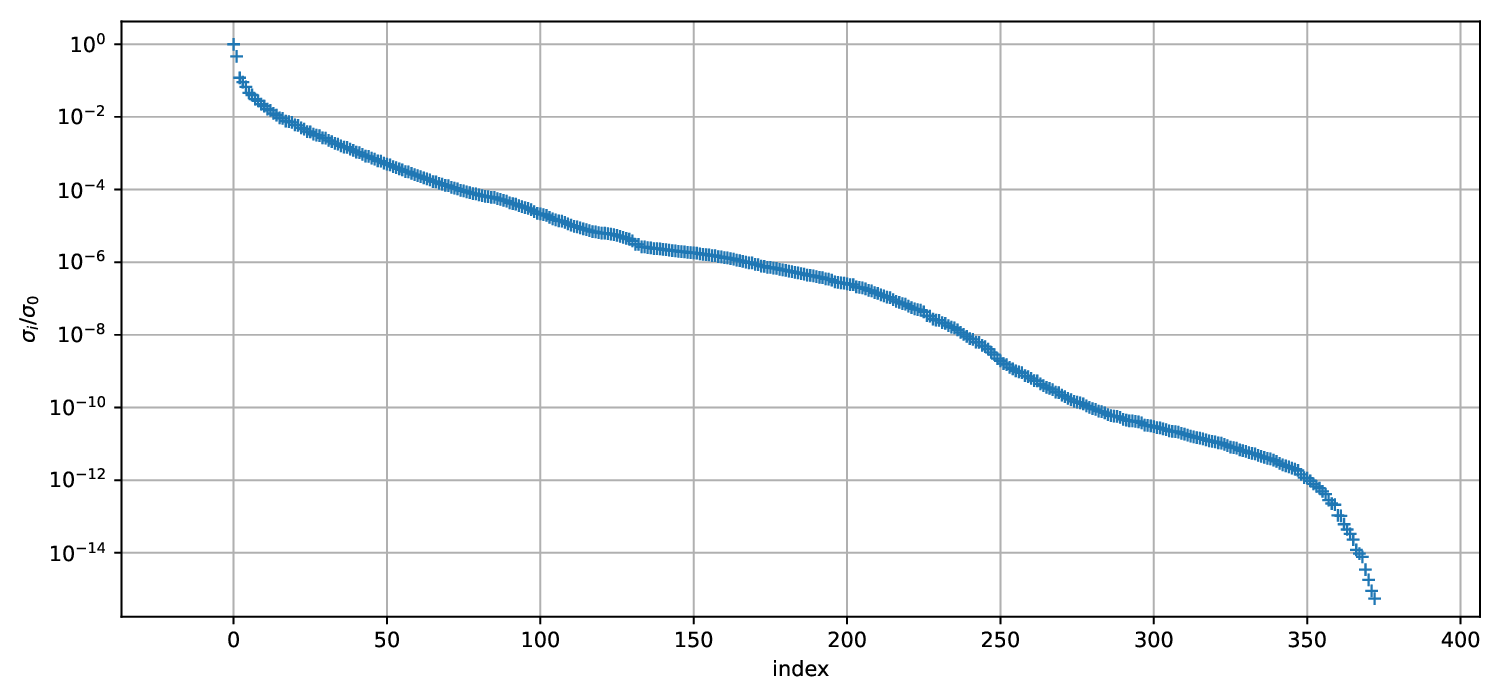}
    \caption{Normalized singular values of the compressible Euler equations data for the modified Sod's shock tube problem.}
    \label{fig:eulersingvals}
\end{figure}

\begin{figure}
    \centering
    \includegraphics[width = \textwidth]{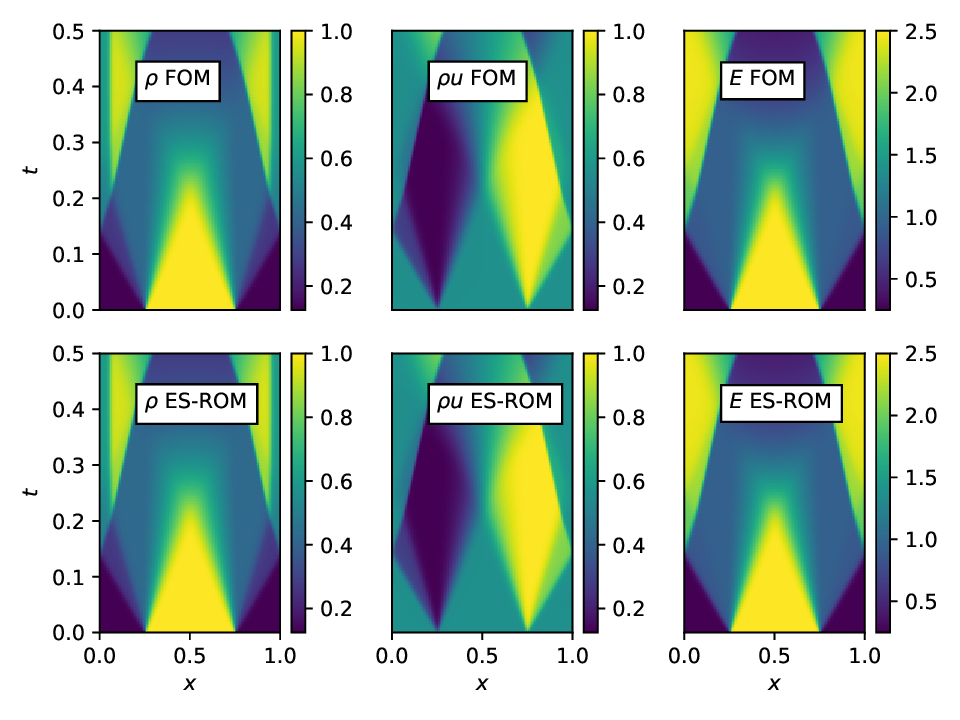}
    \caption{$x-t$ plot of solution approximation by the entropy stable ROM compared to the FOM for the modified Sod's shock tube problem.}
    \label{fig:eulerxt}
\end{figure}

We will be primarily interested in the entropy projection and tangent space enrichment during this experiment, but for completeness we also plot the ROM approximations to the FOM solution and the singular values. The singular values are displayed in \autoref{fig:eulersingvals} and the ROM approximations are shown using $x-t$ plots in \autoref{fig:eulerxt}. A relatively slow decay of singular values can be observed in \autoref{fig:eulersingvals}, hence linear model reduction methods are likely to not perform well for this problem. The FOM solution is approximated well by our novel entropy stable manifold Galerkin ROM. The solution approximations of the ROM as displayed in \autoref{fig:eulerxt} are nearly identical to the FOM. 

\begin{figure}
    \centering
    \includegraphics[width = \textwidth]{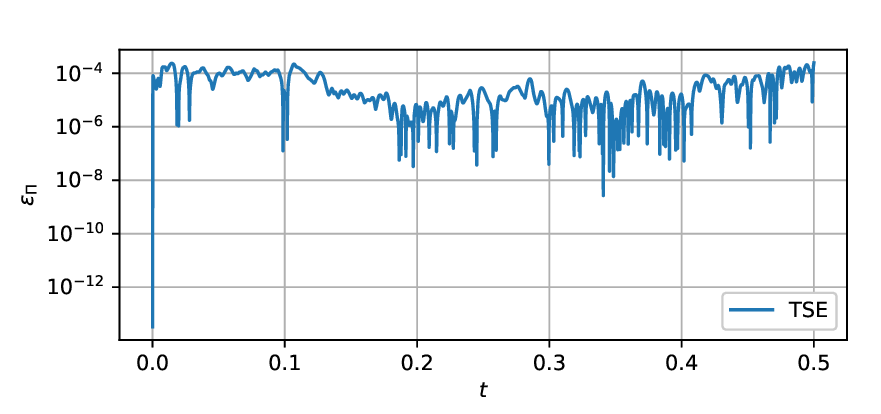}
    \caption{The relative entropy projection error for a ROM with TSE. The ROM without TSE resulted in all NaN values and is therefore not included in the figure.}
    \label{fig:eulerrelentproj}
\end{figure}

\begin{figure}
    \centering
    \includegraphics[width = \textwidth]{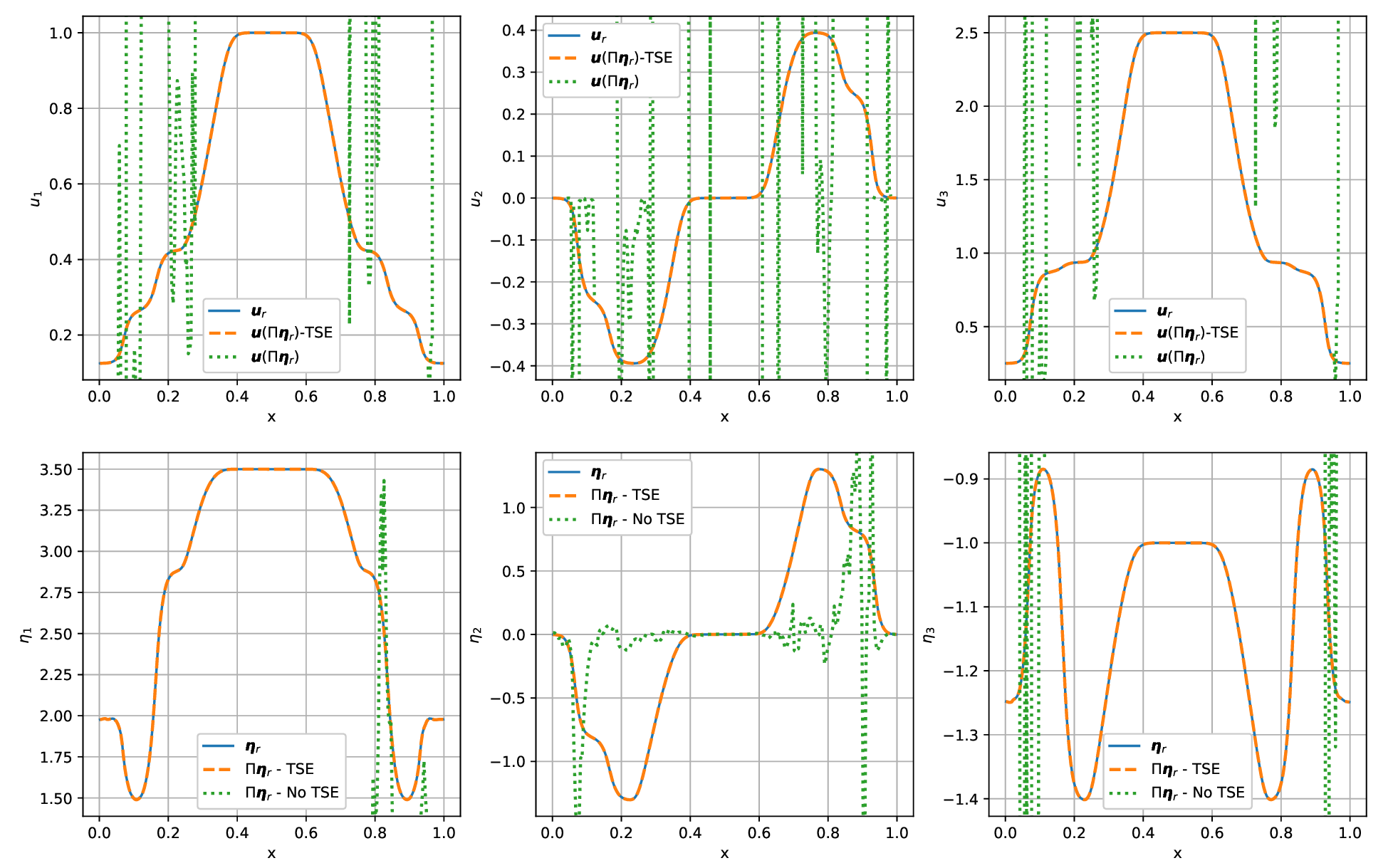}
    \caption{Entropy variable $\bm{\eta}_r$ and conserved variable $\bm{u}_r$ approximation by entropy projection with and without TSE at $t_p = 0.1$.}
    \label{fig:eulerrec01}
\end{figure}

\begin{figure}
    \centering
    \includegraphics[width = \textwidth]{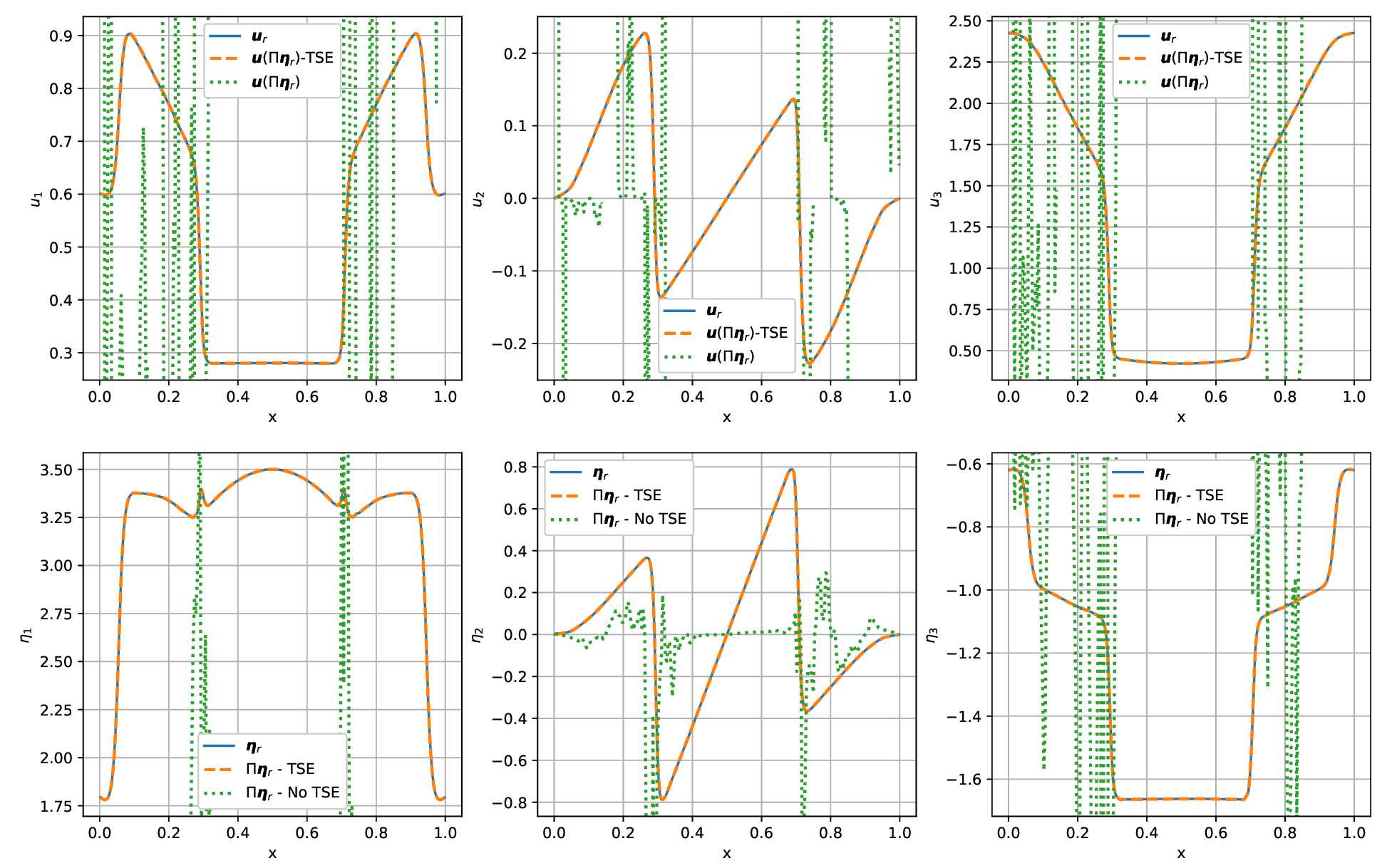}
    \caption{Entropy variable $\bm{\eta}_r$ and conserved variable $\bm{u}_r$ approximation by entropy projection with and without TSE at $t_p = 0.5$.}
    \label{fig:eulerrec05}
\end{figure}

We are interested in the accuracy of the entropy projection with and without tangent space enrichment. Accordingly, we introduce a metric to measure this accuracy. Since we are not only interested in comparing errors, but also to get an idea of the absolute size of the error we specifically introduce the relative entropy projection error:
\begin{equation*}
    \varepsilon_{\Pi}(t) := \frac{||\bm{u}_r(t) - \bm{u}(\Pi_{T\mathcal{M}}\bm{\eta}_r(t))||_{\Omega_h}}{||\bm{u}_r(t)||_{\Omega_h}},
\end{equation*}
measuring not only how far the entropy projection is from the identity mapping as with the entropy projection error $\varepsilon_s$ of \autoref{sec:tse}, but also the size of the error $\varepsilon_s$ relative to the approximated value $\bm{u}_r$. We will plot this value for two ROM simulations of the compressible Euler equations with an entropy projection, where one has an enriched tangent space and the other not. The results are provided in \autoref{fig:eulerrelentproj}. The entropy projection error with TSE can be seen to be a very small fraction of the norm $||\bm{u}_r(t)||_{\Omega_h}$, indicating minimal impact on the accuracy of the ROM given it is well-conditioned. In contrast, the ROM without TSE instantly produces NaN values and could therefore not be included in \autoref{fig:eulerrelentproj}. To obtain a further comparison we plot the spatial profiles of the projected entropy variables at two selected moments $t_p \in \mathbb{R}^+$ in time, namely $t_p \in \{0.1,0.5\}$. To have a meaningful comparison, i.e.\ one where we are not projecting NaN values to start with, we calculate the entropy variables from the stable ROM with enriched tangent space. Furthermore, we use the the generalized coordinates $\bm{a}_p = \Phi^T X_p$ to compute the tangent space basis for the ROM without TSE. The results are shown in \autoref{fig:eulerrec01} and \autoref{fig:eulerrec05}. Very poor reconstruction of entropy variables can be observed for the ROM without TSE, whereas with TSE the reconstruction is accurate at both moments. From \autoref{fig:eulerrec01} and \autoref{fig:eulerrec05} the NaN values in \autoref{fig:eulerrelentproj} can be explained by the projection of the entropy variables taking unphysical values (positive $\eta_3$). From this we conclude that tangent space enrichment or any other manner of assuring the accuracy of the entropy projection is vital for a properly functioning ROM when using an entropy projection.

\begin{figure}
    \centering
    \includegraphics[width = \textwidth]{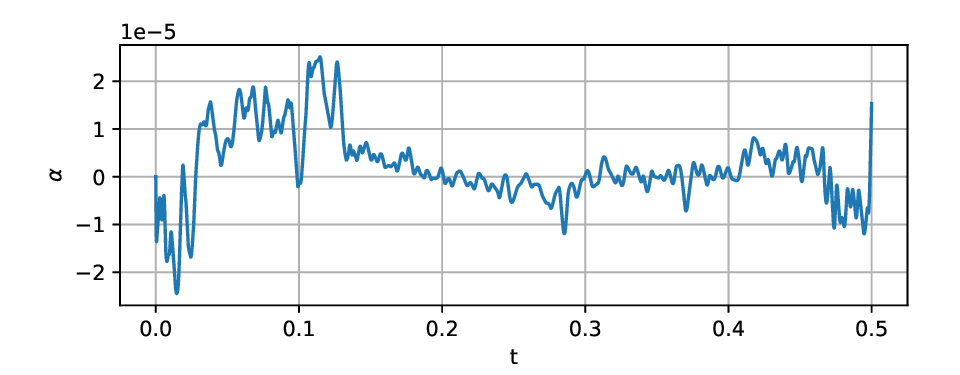}
        \caption{Value of the artificial TSE coordinate $\alpha$ throughout a simulation of the modified Sod's shock tube experiment.}
        \label{fig:euleralpha}
\end{figure}

\begin{figure}
    \centering
    \includegraphics[width = \textwidth]{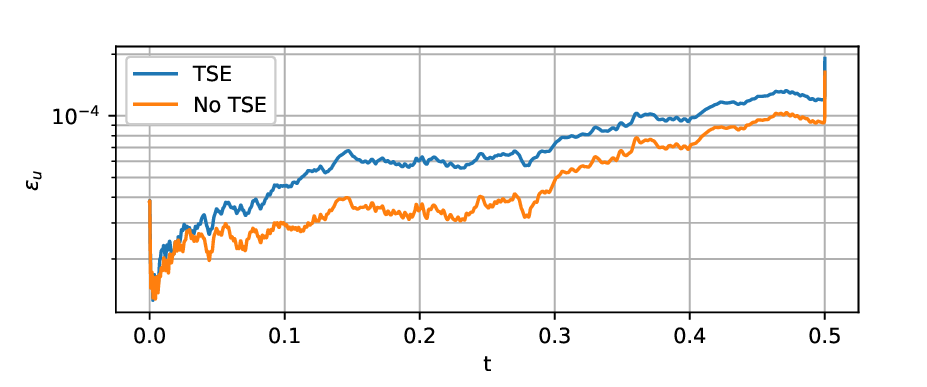}
        \caption{Value of the ROM error $\varepsilon_u$ throughout a simulation of the modified Sod's shock tube experiment.}
        \label{fig:eulerepsu}
\end{figure}

In applying the tangent space enrichment framework, we rely on the artificial TSE coordinate $\alpha$ staying small ($\alpha \ll 1$) during simulations. If this is not the case we cannot assure that the reduced space can accurately represent the solution nor that the local tangent space can accurately represent the FOM dynamical system $\frac{d\bm{u}_h}{dt}$ at that point. The reason for this being that the enriched manifold parameterization $\hat{\bm{\varphi}}$ and its Jacobian matrix $\hat{\bm{J}}$ are no longer close to the original parameterization $\bm{\varphi}$ and Jacobian $\bm{J}$ which are accurate by assumption. To verify this is indeed not the case we will monitor the value of $\alpha$ throughout a simulation of the compressible Euler equations. The results are shown in \autoref{fig:euleralpha}. We have also plotted the error $\varepsilon_u$ of the ROM with entropy projection and tangent space enrichment and of a generic manifold Galerkin ROM for reference in \autoref{fig:eulerepsu}. It can be seen in \autoref{fig:euleralpha} that the value of $\alpha$ remains small around $\mathcal{O}(10^{-5})$. Consequently, the original manifold parameterization $\bm{\varphi}$ and Jacobian $\bm{J}$ are well-approximated by their enriched counterparts $\hat{\bm{\varphi}}$ and $\hat{\bm{J}}$. It can be seen in \autoref{fig:eulerepsu} that this is, in fact, the case, as the errors of the ROM differ by at most $\mathcal{O}(10^{-5})$ and evolve very similarly. Hence, we conclude that the entropy projection step, introduced to obtain an entropy stable framework, is not detrimental to the accuracy of the ROM, provided that an enriched tangent space is used.

\section{Conclusion}\label{sec:conclusion}
In this article we have proposed a method to construct nonlinear manifold Galerkin reduced order models (ROMs) in such a way that the important total entropy functional of the ROM approximation is appropriately conserved or dissipated. This is a crucial concept in obtaining stable and physically admissible ROM solutions. In particular, we have focused on systems of one-dimensional nonlinear conservation laws. Correct semi-discrete entropy estimates upon orthogonal projection were obtained for these systems by evaluating the projected system not at the current state, but at its entropy projection. This was proposed earlier for linear ROMs and extended in this work to nonlinear manifold ROMs. 

The entropy projection of the state is obtained by transforming conservative variables to entropy variables, consequently projecting these on the tangent space of the reduced manifold and finally transforming back to conserved variables. To assure accuracy, it is important that the entropy projection is as close as possible to the identity mapping. This is generally not the case for general nonlinear reduced spaces and hence we have proposed the method of tangent space enrichment (TSE). With TSE the manifold is lifted along an additional dimension parameterized by a new coordinate. This coordinate direction is constructed to linearly extend in the direction of the local entropy variables, so that the tangent space spans the entropy variable at least approximately given the absolute value of the TSE coordinate. Accordingly, the entropy projection error will stay small. 


We have tested our proposed framework on several nonlinear conservation laws from fluid dynamics. We verified that the entropy estimates are satisfied (semi-discretely): the projection of entropy-conserving flux differences produces no total entropy and the projection of entropy dissipative terms dissipates total entropy. This is in contrast to the generic manifold Galerkin framework which leads to production of entropy in our numerical experiments, which is physically incorrect. We have also shown that the introduction of the artificial TSE coordinate  is vital for the accuracy of the entropy projection and leads to minimal decreases in accuracy.

We have also for the first time generalized the recently proposed quadratic manifolds to rational quadratic manifolds. We have suggested a framework to find the coefficients of the rational quadratic manifolds based on a nonlinear curve fitting approach. We have also formulated the rational quadratic polynomials such that they do no not have real poles. This was achieved through semi-definite constraints, avoiding division by zero for any point in the reduced space. Numerical experiments on the inviscid Burgers equation have shown the increased performance of these rational quadratic manifold parameterizations as compared to existing quadratic manifold parameterizations and linear methods.

In future work, two challenges need to be tackled to make the approach computationally efficient: (i) we need a faster way to fit the rational quadratic manifolds and (ii) we need an entropy-stable hyperreduction approach. The former can possibly be achieved through linearization and iterative techniques combined with better choices of generalized coordinates \cite{schwerdtnergreedy}, whereas the latter could be achieved by adapting the constrained optimization formulation that we proposed for energy-conserving systems in \cite{kleinenergy}. In addition, the framework would benefit from extension with an entropy-stable time integration method and entropy-stable treatment of boundary conditions.

\section*{CRediT authorship contribution statement}
\textbf{R.B. Klein}:  conceptualization, methodology, software, validation, formal analysis, investigation, writing - original draft \\
\textbf{B. Sanderse}: writing - review \& editing, supervision \\
\textbf{P. Costa}: writing - review \& editing, supervision \\
\textbf{R. Pecnik}: writing - review \& editing, supervision \\
\textbf{R.A.W.M. Henkes}: writing - review \& editing, supervision, project administration, funding acquisition \\

\section*{Declaration of competing interests}
The authors declare that they have no known competing financial interests or personal relationships that could have appeared to influence the work reported in this paper.

\section*{Data availability}
Data will be made available on request.

\section*{Acknowledgements}
R.B. Klein gratefully acknowledges the funding for this project obtained from Delft University of Technology.


\bibliography{report}{}
\bibliographystyle{plain}

\end{document}